\input amstex
\documentstyle{amsppt}
\NoBlackBoxes

\TagsOnRight

\def\cal{\Cal}
\def\AA{{\cal A}}

\def\HH{{\cal H}}

\def\UU{{\cal U}}
\def\VV{{\cal V}}

\def\QQ{{\cal Q}}

\def\PP{{\cal P}}
\def\Z{{\Bbb Z}}
\def\C{{\Bbb C}}
\def\R{{\Bbb R}}

\def\e{{\epsilon}}

\def\n{\noindent}
\def\part{{\partial}}

\rightheadtext{Hamiltonian homeomorphisms} \leftheadtext{Yong-Geun Oh \& Stefan M\"uller}

\topmatter
\title
The group of Hamiltonian homeomorphisms \\ and $C^0$-symplectic topology
\endtitle
\author
Yong-Geun Oh \footnote{Partially supported by the NSF Grants \#
DMS-0203593 and \# DMS 0503954, Vilas Research Award of University
of Wisconsin and by a grant of the Korean Young Scientist Prize
\hskip8.5cm\hfill} \& Stefan M\"uller
\endauthor
\address
Department of Mathematics, University of Wisconsin, Madison, WI 53706, ~USA \&
Korea Institute for Advanced Study, 207-43 Cheongryangri-dong Dongdaemun-gu Seoul 130-012, KOREA
\endaddress

\abstract The main purpose of this paper is to carry out some of the foundational study
of $C^0${\it -Hamiltonian geometry} and $C^0${\it -symplectic topology}.
We introduce the notion of {\it Hamiltonian topology} on the space of Hamiltonian paths
and on the group of Hamiltonian diffeomorphisms.
We then define the {\it group} $Hameo(M,\omega)$ of {\it Hamiltonian homeomorphisms} such that
$$ Ham(M,\omega) \subsetneq Hameo(M,\omega) \subset Sympeo(M,\omega) , $$
where $Sympeo(M,\omega)$ is the {\it group} of {\it symplectic homeomorphisms}.
We prove that $Hameo(M,\omega)$ is a {\it normal subgroup} of $Sympeo(M,\omega)$ and
contains all the time-one maps of Hamiltonian vector fields of $C^{1,1}$-functions.
We prove that $Hameo(M,\omega)$ is path-connected and so contained in the
identity component $Sympeo_0(M,\omega)$ of $Sympeo(M,\omega)$.

We also prove that the {\it mass flow} of any element from $Hameo(M,\omega)$ vanishes.
In the case of a closed orientable surface, this implies that $Hameo(M,\omega)$ is strictly smaller than the
identity component of the group of area-preserving homeomorphisms
when $M \neq S^2$. For the case of $S^2$, we conjecture that
$Hameo(S^2,\omega)$ is still a proper subgroup of
$Homeo^\Omega_0(S^2) = Sympeo_0(S^2,\omega)$.
\endabstract

\date Revision in March 2006
\enddate


\endtopmatter


\document

\head \bf Contents \endhead

\n 1. Introduction
\smallskip
\n 2. Symplectic homeomorphisms and the mass flow homomorphism
\smallskip
\n 3. Definition of Hamiltonian topology and the Hamiltonian homeomorphism group
\smallskip
\n 4. Basic properties of the group of Hamiltonian homeomorphisms
\smallskip
\n 5. The two dimensional case
\smallskip
\n 6. The non-compact case and open problems
\smallskip
\n Appendix

\bigskip

\head  \bf \S1. Introduction
\endhead

Let $(M,\omega)$ be a connected symplectic manifold.
{\it Unless explicit mention is made to the contrary, $M$ will be closed.}
See section 6 for the necessary changes in the non-compact case or in the case with boundary.
Denote by $Symp(M,\omega)$ the group of symplectic diffeomorphisms, i.e.,
the subgroup of $Diff(M)$ consisting of diffeomorphisms $\phi:M
\to M$ such that $\phi^*\omega = \omega$. We provide the
$C^\infty$-topology on $Diff(M)$ under which $Symp(M,\omega)$
forms a closed topological subgroup. We call the induced topology
on $Symp(M,\omega)$ the $C^\infty$-topology of $Symp(M,\omega)$. We denote
by $Symp_0(M,\omega)$ the path-connected component of the identity in $Symp(M,\omega)$.
The celebrated $C^0$-rigidity theorem by Eliashberg [El], [Gr] in
symplectic topology states

\proclaim{[$C^0$-Symplectic Rigidity, El]} The subgroup
$Symp(M,\omega) \subset Diff(M)$ is closed in the $C^0$-topology.
\endproclaim

Therefore it is reasonable to define a {\it symplectic
homeomorphism} as any element from
$$ \overline{Symp(M,\omega)} \subset Homeo(M) , $$
where the closure is taken inside the group $Homeo(M)$ of
homeomorphisms of $M$ with respect to the $C^0$-topology
(or compact-open topology).
This closure forms a group and is a
topological group with respect to the induced $C^0$-topology. We
refer to section 2 for the precise definition of the $C^0$-topology on
$Homeo(M)$.

\definition{Definition 1.1 [Symplectic homeomorphism group]}
We denote the above closure equipped with the $C^0$-topology by
$$ Sympeo(M,\omega):= \overline{Symp(M,\omega)} , $$
and call this group the {\it symplectic homeomorphism group}.
\enddefinition

We provide two justifications for this definition.

Firstly, it is easy to see that any symplectic homeomorphism
preserves the {\it Liouville measure} induced by the volume form
$$
\Omega = \frac{1}{n!}\omega^n,
$$
which is an easy consequence of Fatou's lemma in measure theory.
In fact, this measure-preserving property follows from a general
fact that the set of measure-preserving homeomorphisms is closed
in the group of homeomorphisms under the compact-open topology. In
particular in two dimensions, $Sympeo(M,\omega)$ coincides with
$Homeo^\Omega(M)$, where $Homeo^\Omega(M)$ is the group of
homeomorphisms that preserve the Liouville measure. This follows
from the fact that any area-preserving homeomorphism can be
$C^0$-approximated by an area-preserving diffeomorphism in two
dimensions (see Theorem 5.1). Secondly, it is easy to see from
Eliashberg's rigidity that we have
$$
Sympeo(M,\omega) \subsetneq Homeo^\Omega(M) \tag 1.1
$$
when $\dim M \geq 4$.
In this sense the symplectic homeomorphism group is a good high
dimensional {\it symplectic} generalization of the group of area-preserving homeomorphisms.

There is another smaller subgroup $Ham(M,\omega) \subset
Symp_0(M,\omega)$,  the {\it Hamiltonian diffeomorphism group},
which plays a prominent role in many problems in the development
of symplectic topology, starting implicitly from Hamiltonian
mechanics and more conspicuously from the Arnold conjecture. One
of the purposes of the present paper is to give a precise
definition of the $C^0$-counterpart of $Ham(M,\omega)$. This
requires some lengthy discussion on the Hofer geometry of
Hamiltonian diffeomorphisms.

The remarkable Hofer norm of Hamiltonian diffeomorphisms
introduced in [H1,2] is defined by
$$
\|\phi\| = \inf_{H\mapsto \phi} \|H\|, \tag 1.2
$$
where $H \mapsto \phi$ means that $\phi= \phi_H^1$ is the time-one
map of Hamilton's equation
$$ \dot x = X_H(t,x) . $$
In other words, the family $\phi_H^t$ of diffeomorphisms of $M$ satisfies
$$ \frac{d}{dt} \phi_H^t = X_H \circ \phi_H^t, \quad \phi_H^0 = id, $$
i.e., $(t,x) \mapsto \phi_H^t (x)$ is the flow of the Hamiltonian vector field $X_H$
associated to the Hamiltonian function $H : [0,1] \times M \to \R$, defined by
$ X_H \rfloor \omega = dH$,
and $\phi$ is the time-1 map of this flow.
The norm $\|H\|$ is  defined by
$$
\|H\| = \int _0^1 \text{osc }H_t \, dt
= \int_0^1 \left( \max_{x \in M} H_t (x) - \min_{x \in M} H_t (x) \right) \, dt. \tag 1.3
$$
This is a version of the
$L^{(1,\infty)}$-norm on $C^\infty([0,1] \times M,\R)$.

Here $(M,\omega)$ is a general symplectic manifold, which may be
open or closed.
We will always assume that $X_H$ is compactly supported in $Int(M)$ when $M$ is
open so that the flow exists for all time and is supported in
$Int(M)$. For the closed case, we will always assume that the
Hamiltonians are normalized by
$$
\int_M H_t \, d\mu = 0, \quad \text{for all } t \in [0,1],
$$
where $d\mu$ is the Liouville measure. We call such Hamiltonian
functions {\it normalized}.
In both cases, there is a one-one correspondence between $H$ and the path
$\phi_H \colon t \mapsto \phi_H^t$.
There is the $L^\infty$-version of the
Hofer norm originally adopted by Hofer [H1] and defined by
$$
\|H\|_\infty: = \max_{(t,x)} H(t,x) - \min_{(t,x)} H(t,x).
$$
Although this $L^\infty$-norm would be easier to handle and enough
for most of the geometric purposes in the smooth category, we would
like to emphasize that it is important to use the
$L^{(1,\infty)}$-norm (1.3) for the purpose of working with the $C^0$-category:
One essential point that distinguishes the
$L^{(1,\infty)}$-norm from the $L^\infty$-norm is that the
important {\it boundary flattening procedure} is
$L^{(1,\infty)}$-continuous but {\it not} $L^\infty$-continuous.
(See section 3 and Appendix 2 for more precise remarks.)
Recall that this flattening procedure is crucial for defining the
Floer homology and so the spectral invariants [Oh4] and for the
various constructions involving concatenation in symplectic
geometry. Because of this, we adopt the $L^{(1,\infty)}$-norm in
our exposition from the beginning.

When we do not explicitly mention otherwise, we always assume that all the
functions and diffeomorphisms are smooth. In particular,
$Ham(M,\omega)$ is a subgroup of $Symp_0(M,\omega)$.
Banyaga [Ba] proved that this group is a simple group. Recently Ono
[On] gave a proof of the $C^\infty$-Flux Conjecture which implies
that $Ham(M,\omega)$ is a closed subgroup of $Symp_0(M,\omega)$
and locally contractible in the $C^\infty$-topology. The question
whether $Ham(M,\omega)$ is $C^0$-closed in $Symp_0(M,\omega)$ is
sometimes called the $C^0$-Flux Conjecture.

The above norm $\| H \|$ can be identified with the Finsler length
$$
\text{leng}(\phi_H) =
\int_0^1 \Big(\max_{x \in M} H(t,(\phi_H^t)(x)) - \min_{x \in M} H(t,(\phi_H^t)(x))\Big)\, dt \tag 1.4
$$
of the path $\phi_H:t \mapsto \phi_H^t$ where the Banach norm on
$T_{id}Ham(M,\omega) \cong C^\infty(M)/\R$ is defined by
$$
\|h\| = \text{osc}(h) = \max h - \min h
$$
for a normalized function $h:M \to \R$.

\definition{Definition 1.2} We call a continuous
path $\lambda:[0,1] \to Symp(M,\omega)$ a (smooth) {\it Hamiltonian path}
if it is generated by the flow of $\dot x = X_H(t,x)$
with respect to a smooth Hamiltonian $H :[0,1] \times M \to \R$ (see also Definition A.1).
We denote by $\PP^{ham}(Symp(M,\omega))$ the set of
Hamiltonian paths $\lambda$ and by $\PP^{ham}(Symp(M,\omega),id)$
the set of Hamiltonian paths $\lambda$ that satisfy $\lambda(0) = id$.
We also denote by
$$ ev_1:\PP^{ham}(Symp(M,\omega),id) \to Symp(M,\omega) \tag 1.5 $$
the evaluation map $ev_1(\lambda) = \lambda(1) = \phi_H^1$.
\enddefinition
For readers' convenience, we will give a precise description of the
$C^\infty$-topology on $\PP^{ham}(Symp(M,\omega),id)$ in Appendix 1.
By definition, $Ham(M,\omega)$ is the set of images of $ev_1$. We
will be mainly interested in the Hamiltonian paths lying in the
identity component $Symp_0(M,\omega)$ of $Symp(M,\omega)$.

\definition{Definition 1.3 [The Hofer topology]} Consider the metric
$$d_H: \PP^{ham}(Symp(M,\omega),id) \to \R_{\ge 0}$$
defined by
$$ d_H(\lambda,\mu) := \text{leng}(\lambda^{-1}\circ \mu) , \tag 1.6 $$
where $\lambda^{-1}\circ \mu$ is the Hamiltonian path $t \in [0,1]
\mapsto \lambda(t)^{-1}\mu(t)$. We call the induced topology on
$\PP^{ham}(Symp(M,\omega),id)$ the {\it Hofer topology}. We define the Hofer
topology on $Ham(M,\omega)$ to be the strongest topology for which
the evaluation map (1.5) is continuous.
\enddefinition
It is easy to see that this definition of the Hofer topology on
$Ham(M,\omega)$ coincides with the usual one induced by (1.2),
which also shows that the Hofer topology is metrizable.
Of course nontriviality of the topology is not a trivial matter which was
proven by Hofer [H1] for $\C^n$, by Polterovich [P1] for rational
symplectic manifolds and by Lalonde and McDuff in its
complete generality [LM].
It is also immediate to check that the Hofer topology is locally path-connected.

The relation between the Hofer topology on $Ham(M,\omega)$ and the
$C^\infty$-topology or the $C^0$-topology thereon is rather delicate.
However it is known (see [P2] and Example 4.2) that the Hofer norm function
$$
\phi \in Ham(M,\omega) \to \|\phi\|
$$
is {\it not}  continuous with respect to the $C^0$-topology in general.
We refer to [Si], [H2]
for some results for compactly supported Hamiltonian
diffeomorphisms on $\R^{2n}$ in this direction.

The main purpose of this paper is to carry out a foundational
study of $C^0$-Hamiltonian geometry. We first give the precise
definition of a topology on the space of Hamiltonian paths with
respect to which the spectral invariants for Hamiltonian paths
constructed in [Oh3-6] will all be continuous [Oh7]. We then
define the notion of {\it Hamiltonian homeomorphisms} and denote
the set thereof by $Hameo(M,\omega)$. We provide many evidences
for our thesis that the Hamiltonian topology is the right topology
for the study of topological Hamiltonian geometry. In fact, the
notion of Hamiltonian topology has been vaguely present in the
literature without much emphasis on its significance (see [H2],
[V], [HZ], [Oh3] for some theorems related to this topology). {\it
However all of the previous works fell short of constructing a
``group'' of continuous Hamiltonian maps}. A precise formulation
of the topology will be essential in our study of the continuity
property of spectral invariants, and also in our construction of
$C^0$-symplectic analogs corresponding to various
$C^\infty$-objects or invariants. We refer readers to [Oh7] for
the details of this study.

The following is the $C^0$-analog to the well-known fact that
$Ham(M,\omega)$ is a normal subgroup of $Symp_0(M,\omega)$.
\proclaim{Theorem I} The group $Hameo(M,\omega)$ forms a normal subgroup of $Sympeo(M,\omega)$.
\endproclaim
We also prove
\proclaim{Theorem II} $Hameo(M,\omega)$ is path-connected and
contained in the identity component of $Sympeo(M,\omega)$, i.e., we have
$$ Hameo(M,\omega) \subset Sympeo_0(M,\omega). $$
\endproclaim
See Theorems 4.4 and 4.5 respectively.
In section 4, we also prove that all Hamiltonian diffeomorphisms
generated by $C^{1,1}$-Hamiltonian functions are contained in
$Hameo(M,\omega)$ and give an example of a Hamiltonian homeomorphism
that is not even Lipschitz (see Theorem 4.1 and Example 4.2 respectively).
We recall the notion of the {\it mass flow homomorphism} [S], [T], [Fa],
which is also called {\it the mean rotation vector} in the
literature on area-preserving maps.

We prove (see Theorem 5.2 and Theorem 5.5.)

\proclaim{Theorem III}  The values of
the mass flow homomorphism with respect to the Liouville
measure of $\omega$ are zero on $Hameo(M,\omega)$.
\endproclaim
As a corollary to Theorems I - III, we prove that in dimension two $Hameo(M,\omega)$
is strictly smaller than the identity component of the group of
area-preserving homeomorphisms if $M \neq S^2$. For the case of
$S^2$, we still conjecture
\proclaim{Conjecture 1} Let $M=S^2$
with the standard area form $\omega = \Omega$. $Hameo(S^2,\omega)$ is a
proper subgroup of $Homeo^\Omega_0(S^2) = Sympeo_0(S^2,\omega)$.
\endproclaim

The last equality follows from Theorem 5.1.
Therefore one consequence of Conjecture 1 together with normality (Theorem I) and
path-connectedness (Theorem II) would be the following result,
which would answer negatively to the following open question
since the work of Fathi [Fa] appeared.
\proclaim{Conjecture 2}
$Homeo^\Omega_0(S^2)$, the identity component of the group of area-preserving
homeomorphisms of $S^2$, is not a simple group.
\endproclaim

We refer to section 5 for further discussions on the relation
between $Hameo(M,\omega)$ and the simpleness question of the area-preserving homeomorphism group of $S^2$.

In section 6, we look at the open case and define the
corresponding Hamiltonian topology and the $C^0$-version of compactly
supported Hamiltonian diffeomorphisms.

Finally we have two appendices. In Appendix 1, we provide precise
descriptions of the $C^\infty$-topologies on $Ham(M,\omega)$ and
its path space $\PP^{ham}(Symp(M,\omega),id)$. We also give the
proof of the fact that $C^\infty$-continuity of a Hamiltonian path
implies the continuity with respect to the Hamiltonian topology.
In Appendix 2, we recall the proof of the
$L^{(1,\infty)}$-Approximation Lemma from [Oh3] in a more precise
form for the readers' convenience.

The senior author is greatly indebted to the graduate students of
Madison attending his symplectic geometry course in the fall of
2003. He thanks them for their patience listening to his
lectures throughout the semester, which were sometimes erratic in
some foundational materials concerning the Hamiltonian
diffeomorphism group. The present paper partly grew out of the
course. He also thanks J. Franks, J. Mather and A. Fathi for a
useful communication concerning the smoothing of area-preserving
homeomorphisms. Writing of the original version of this paper has
been carried out while the senior author was visiting the Korea
Institute for Advanced Study in the winter of 2003. He thanks KIAS
for its financial support and excellent research atmosphere.

We thank A. Fathi for making numerous helpful comments on a
previous senior author's version of the paper, which has led to
corrections of many erroneous statements and proofs and to
streamlining the presentation of the paper. We also thank the
referee for carefully reading the previous version and pointing out
many inaccuracies, and for providing many helpful suggestions on improving the
presentation of the paper.

During the preparation of the current revision, Viterbo [V2]
answered affirmatively to the $C^0$-version of Question 3.16,
and subsequently the senior author proved its $L^{(1,\infty)}$-version [Oh7].

\bigskip

\n{\bf Notations} \medskip \roster
\item Unless otherwise stated,
$H$ always denotes a {\it normalized smooth} Hamiltonian function $[0,1] \times M \to \R$, and
we always denote by $\|\cdot\|$ the $L^{(1,\infty)}$-norm
$$
\|H\| = \int_0^1 \left(\max_{x \in M} H (t,x) - \min_{x \in M} H (t,x)\right) \,dt .
$$
We denote by $C_m^\infty ([0,1] \times M,\R)$ the space of such functions $H$ with the norm $\|\cdot\|$,
and by $L_m^{(1,\infty)} ([0,1] \times M,\R)$ its completion with respect to $\|\cdot\|$.
\item Our convention is that $\phi_H$ always denotes a smooth Hamiltonian path $\phi_H : t \mapsto \phi_H^t$,
while $\phi$ or $\phi_H^t$ denotes a single diffeomorphism.
Unless otherwise stated, $\| \phi \|$ always denotes the Hofer
norm (1.2) for $\phi \in Ham(M,\omega)$.
\item $G_0$: The identity path-component of any topological group $G$.
\item $Homeo(M)$: The group of homeomorphisms of $M$ with the $C^0$-topology.
We will often abbreviate composition of maps by $\psi \circ \phi = \psi \phi$.
\item $\PP(G)$, \, $\PP(G,id)$: The space of continuous paths $\lambda : [0,1] \to G$,
and the space of continuous paths with $\lambda (0) = id$, respectively.
\item $Homeo^\Omega(M)$: The topological subgroup of $Homeo(M)$
consisting of measure (induced by the volume form $\Omega$)
preserving homeomorphisms of $M$.
\item $Symp(M,\omega)$: The group
of symplectic diffeomorphisms with the $C^\infty$-topology.
\item $Sympeo(M,\omega)$: The $C^0$-closure of $Symp(M,\omega)$ in $Homeo(M)$.
\item $\PP^{ham}(Symp(M,\omega),id)$: The space of smooth Hamiltonian paths
$\lambda:[0,1]\to Symp(M,\omega)$ with $\lambda(0) = id$ with the $C^\infty$-topology.
\item $\PP^{ham}_s (Symp(M,\omega),id)$: $\PP^{ham}(Symp(M,\omega),id)$ with the (strong) Hamiltonian topology.
\item $Ham(M,\omega) \subset Symp_0(M,\omega)$: The subgroup of Hamiltonian
diffeomorphisms with the $C^\infty$-topology.
\item $\HH am(M,\omega)$: $Ham(M,\omega)$ with the (strong) Hamiltonian topology.
\item $ev_1 : \PP^{ham}(Symp(M,\omega),id) \to Ham(M,\omega)$ the evaluation map.
\item $Hameo(M,\omega)$: The group of (strong) Hamiltonian homeomorphisms with
the $C^0$-topology.
\item $\HH ameo(M,\omega)$: $Hameo(M,\omega)$ with the (strong) Hamiltonian topology.
\endroster

\head{\bf \S 2. Symplectic homeomorphisms and the mass flow homomorphism}
\endhead

Let $(M,\omega)$ be as in the introduction.
We fix any Riemannian metric and denote by $d$ the induced
Riemannian distance function on $M$.
We denote by $Homeo_0(M)$ the path-connected component of the identity in $Homeo(M)$,
the group of homeomorphisms of $M$.
Denote by $\PP(Homeo(M),id)$ the
set of continuous paths $\lambda:[0,1] \to Homeo(M)$ with $\lambda(0)=id$.
We denote by $d_{C^0}$ the standard $C^0$-distance of {\it maps} defined by
$$
d_{C^0}(\phi,\psi) = \max_{x\in M}\left(d(\phi(x),\psi(x)\right).
$$
Then for any two {\it homeomorphisms} $\phi, \, \psi \in Homeo(M)$ we define their $C^0$-distance
$$
\overline d(\phi,\psi) =
\max \left \{d_{C^0}(\phi,\psi),d_{C^0}(\phi^{-1},\psi^{-1}) \right \} . \tag 2.1
$$
With respect to this metric, $Homeo(M)$ becomes a complete metric space.
We call the topology induced by $\overline{d}$ the $C^0${\it -topology} on $Homeo(M)$.
It is easy to see that this topology coincides with the compact-open topology.
In particular, it does not depend on the choice of the particular Riemannian metric.
As we defined in Definition 1.1 of the introduction, the
symplectic homeomorphism group $Sympeo(M,\omega)$ is defined to be
the closure of $Symp(M,\omega)$ in $Homeo(M)$ with respect to this
metric.

Then for given continuous paths $\lambda, \, \mu: [0,1] \to Homeo_0(M)$
with $\lambda(0) = \mu(0)=id$, we define their $C^0$-distance by
$$
\overline d(\lambda,\mu): = \max_{t \in
[0,1]}\overline d(\lambda(t),\mu(t)), \tag 2.2
$$
and call the induced metric topology the $C^0${\it -topology} on $\PP(Homeo(M),id)$.

If $\psi_i$ is a Cauchy sequence in the $C^0$-topology
converging to a homeomorphism $\psi \in Homeo(M)$,
we will write $\lim_{C^0} \psi_i = \psi$.
It is easy to see that $\lim_{C^0} \psi_i^{-1} = \psi^{-1}$
and $\lim_{C^0} \psi_i \phi_i = \psi \phi$ for two sequences
$\lim_{C^0} \psi_i = \psi$ and $\lim_{C^0} \phi_i = \phi$.
The same observations hold for the complete metric (2.2) for continuous paths.
More precisely, let $\lambda_i$ and $\mu_i \in \PP (Homeo(M),id)$
be two Cauchy sequences of continuous paths.
Then there exist {\it continuous} paths $\lambda = \lim_{C^0} \lambda_i \in \PP (Homeo(M),id)$,
$\mu = \lim_{C^0} \mu_i \in \PP (Homeo(M),id)$,
and we have $\lim_{C^0} \lambda_i \mu_i = \lambda \mu$ and
$\lim_{C^0} \lambda_i^{-1} = \lambda^{-1}$.
Here $\lambda^{-1} : [0,1] \to Homeo(M)$ denotes the path $t \mapsto (\lambda(t))^{-1}$.
We will use this frequently in sections 3 and 4.

\medskip

Recall that the symplectic form $\omega$ induces a measure on $M$
by integrating the volume form
$$
\Omega = \frac{1}{n!}\omega^n.
$$
We will call the induced measure the {\it Liouville measure} on
$M$. We denote the Liouville measure by $d\mu =d\mu^\omega$.

The following is an immediate consequence of the well-known fact
(see [Corollary 1.6, Fa] for example) that for any given finite
Borel measure $d\mu$, the group of measure-preserving
homeomorphisms is closed under the above compact-open topology.

\proclaim{Proposition 2.1} Any symplectic homeomorphism $h \in
Sympeo(M,\omega)$ preserves the Liouville measure. More precisely,
$Sympeo(M,\omega)$ forms a closed subgroup of $Homeo^\Omega(M)$.
\endproclaim

It is easy to derive from Eliashberg's rigidity theorem
the properness of the subgroup $Sympeo(M,\omega)
\subset Homeo^\Omega(M)$ when $\dim M \geq 4$.

Next we briefly review the construction from [Fa] of the {\it mass
flow homomorphism} for measure-preserving homeomorphism. When
considered on an orientable surface, it coincides with the
symplectic flux (up to Poincar\'e duality), and it will be used in
section 5 to prove, when $M \neq S^2$, that $Sympeo_0(M,\omega)$
is strictly bigger than the group $Hameo(M,\omega)$ of Hamiltonian
homeomorphisms which we will introduce in the next section.

Let $\Omega$ be a volume form on $M$ and denote by $ Homeo^\Omega_0(M)$
the path-connected component of the identity in the set of measure (induced by $\Omega$) preserving
homeomorphisms with respect to the $C^0$-topology (or compact-open topology).
By Proposition 2.1, we have the inclusion $Sympeo(M,\omega) \subset Homeo^\Omega(M)$.
We will not be studying this inclusion carefully here except in two dimensions.

For any $G$ one of the above groups, we will denote by $\PP(G)$
(respectively  $\PP(G,id)$) the space of continuous path from
$[0,1]$ into $G$ (respectively with $c(0) = id$) with the induced
$C^0$-topology. We denote by $c=(h_t):[0,1] \to G$ the
corresponding path. Since $Homeo^\Omega(M)$ is locally
contractible [Fa], the universal covering space of
$Homeo_0^\Omega(M)$ is represented by homotopy classes of paths $c
\in \PP(Homeo_0^\Omega(M),id)$ with fixed end points. We denote by
$$
\pi : \widetilde{Homeo_0^\Omega}(M) \to Homeo_0^\Omega(M)
$$
the universal covering space and by $[c]$ the corresponding
elements. To define the mass flow homomorphism
$$
\widetilde \theta: \widetilde{Homeo_0^\Omega}(M) \to H_1(M,\R),
\tag 2.3
$$
we use the fact that $H_1(M,\R) \cong \text{Hom}([M,S^1],\R)$,
where $[M,S^1]$ is the set of homotopy classes of maps from $M$ to
$S^1$.

Denote by $C^0(M,S^1)$ the set of continuous maps $M \to S^1$
equipped with the $C^0$-topology. Note that $C^0(M,S^1)$ naturally
forms a group. Identifying $S^1$ with $\R/\Z$, write the group law
on $S^1$ additively. Given $c=(h_t) \in
\PP(Homeo_0^\Omega(M),id)$, we define a continuous group
homomorphism
$$
\widetilde\theta(c): C^0(M,S^1) \to \R
$$
in the following way: let $f: M \to S^1 = \R/\Z$ be continuous.
The homotopy $fh_t - f: M\to S^1$ satisfies $fh_0 - f = 0$, hence
we can lift it to a homotopy $\overline{fh_t - f}: M \to \R$ such
that $\overline{fh_0 - f} = 0$. Then we define
$$
\widetilde\theta(c)(f) = \int_M \overline{fh_1 -f} \, d\mu ,
$$
where $d\mu$ is the given measure on $M$. This induces a
homomorphism
$$
\widetilde \theta : \PP(Homeo_0^\Omega(M),id) \to
Hom(C^0(M,S^1),\R). \tag 2.4
$$
One can check that for each given $f \in C^0(M,S^1)$, the
assignment $c \mapsto \widetilde\theta(c)(f)$ is continuous, i.e.,
the map (2.4) is {\it weakly continuous}. Furthermore
$\widetilde\theta(c)(f)$ depends only on the homotopy class of
$f$, $\widetilde\theta(c)$ is a homomorphism,
$\widetilde\theta(c)$ depends only on the equivalence class of
$c$, and $\widetilde\theta$ is a homomorphism [Fa]. Therefore it
induces a group homomorphism (2.3). The weak continuity of (2.4)
then induces the continuity of the map (2.3).

If we put
$$
\Gamma = \widetilde\theta\left(\ker \left( \pi \colon
\widetilde{Homeo_0^\Omega}(M) \to
Homeo_0^\Omega(M)\right) \right) ,
$$
we obtain by passing to the quotient a group homomorphism
$$
\theta:Homeo^\Omega_0(M) \to H_1(M,\R)/\Gamma, \tag  2.5
$$
which is also called the {\it mass flow homomorphism}.
The group $\Gamma$ is shown to be discrete because it is contained
in $H_1(M,\Z)$ (after normalizing $\Omega$ so that $\int_M \Omega = 1$)
[Proposition 5.1, Fa].

We summarize the above discussion and some fundamental results of
Fathi [Fa] restricted to the case where $M$ is a (smooth)
manifold. Note that Fathi equips $\PP(Homeo(M),id)$ with the
compact-open topology, while we use the $C^0$-topology (2.2). It
is easy to see that the $C^0$-topology is stronger than the
compact-open topology on the path space $\PP(Homeo(M),id)$, and
therefore Fathi's results also apply to our case.

\proclaim{Theorem 2.2 [Fa]} Suppose that $M$ is a closed smooth
manifold and $\Omega$ is a volume form on $M$. Then \roster \item
$Homeo^\Omega(M)$ is locally contractible, \item the map
$\widetilde \theta$ in (2.4) is weakly continuous, and $\theta$ in
(2.5) is continuous, with respect to the $C^0$-topology,
\item the
map $\widetilde \theta$ in (2.3) is surjective, and hence so is
$\theta$, \item $\ker \theta = [\ker \theta, \ker \theta]$ is
perfect, and $\ker \theta$ is simple, if $n \geq 3$.
\endroster
\endproclaim

The following still remains an open problem concerning the
structure of the area-preserving homeomorphism groups in two
dimensions (note that since $H_1(S^2, \R) = 0$, we have $\ker
\theta = Homeo_0^{\Omega}(S^2)$)

\definition{Question 2.3} Is $\ker \theta$ simple when $n = 2$?
In particular, is $Homeo^\Omega_0(S^2)$ a simple group?
\enddefinition

\head{\bf \S 3. Definition of Hamiltonian topology and the Hamiltonian homeomorphism group}
\endhead

We start by recalling the following proposition proven by the
senior author [Oh3] in relation to his study of the length minimizing
property of geodesics in Hofer's Finsler geometry on
$Ham(M,\omega)$. This result was the starting point of the senior
author's research carried out in this paper.

\proclaim{Proposition 3.1 [Lemma 5.1, Oh3]}
Let $\phi_{G_i}$ be a sequence of smooth Hamiltonian paths and
$\phi_G$ be another smooth Hamiltonian path such that
\roster
\item each $\phi_{G_i}$ is length minimizing in its homotopy class
relative to the end points,
\item $\text{\rm leng}(\phi_G^{-1} \phi_{G_i}) \to 0$ as
$i\to \infty$, and
\item the sequence of Hamiltonian paths $\phi_{G_i}$ converges to
$\phi_G$ in the $C^0$-topology.
\endroster
Then $\phi_G$ is length minimizing in its homotopy class
relative to the end points.
\endproclaim

In fact, an examination of the proof of Lemma 5.1 in [Oh3]
shows that the same holds even without (3).
This proposition can be translated into the statement that
the length minimizing property of Hamiltonian paths in its homotopy class
{\it relative to the end points} is closed under a certain
topology on the space of Hamiltonian paths.
In this section, we will first introduce the corresponding topology
on the space of Hamiltonian paths. Then using this topology, which we
call {\it Hamiltonian topology}, we will construct
the group of {\it Hamiltonian homeomorphisms}.

We first recall the definition of ($C^\infty$-)Hamiltonian
diffeomorphisms (see also section 1): A
$C^\infty$-diffeomorphism $\phi$ of $(M,\omega)$ is {\it
$C^\infty$-Hamiltonian} if $\phi = \phi_H^1$ for a
$C^\infty$-function $H: [0,1] \times M \to \R$. Here $\phi_H^1$ is again
the time-one map of the Hamilton equation
$$
\dot x = X_H(t,x).
$$
We denote the set of Hamiltonian diffeomorphisms by $Ham(M,\omega)$,
and recall that $Ham(M,\omega) \subset Symp_0(M,\omega)$.
We will always denote by $\phi_H$ the corresponding Hamiltonian path
$\phi_H: t \mapsto \phi_H^t$ generated by the Hamiltonian $H$
and by $H \mapsto \phi$ when $\phi = \phi_H^1$. In the latter
case, we also say that the diffeomorphism $\phi$ is generated by the
Hamiltonian $H$.

We recall that for two Hamiltonian functions $H$ and $K$,
the Hamiltonian $H\# K$ is given by the formula
$$
(H \# K)_t = H_t + K_t \circ (\phi_H^t)^{-1} \tag 3.1
$$
and generates the path $\phi_{H} \phi_K : t\mapsto \phi_H^t\phi_K^t$.
And the inverse Hamiltonian $\overline H$ corresponding to the
inverse path $(\phi_H)^{-1} : t \mapsto (\phi_H^t)^{-1}$ is defined by
$$ (\overline H)_t = - H_t \circ \phi_H^t . \tag 3.2 $$
We also recall that the Hamiltonian $\psi^* H$,
$$
(\psi^*H)_t = H_t \circ \psi, \tag 3.3
$$
generates the path $\psi^{-1} \phi_H \psi : t \mapsto
\psi^{-1} \phi_H^t \psi$ for any $\psi \in Symp(M,\omega)$. In
particular, $Ham(M,\omega)$ is a normal subgroup of
$Symp(M,\omega)$.
We will be mainly interested in paths of the form
$\phi_H^{-1}\phi_K$. By the above, this path is generated by
$\overline{H} \# K$, and
$$ (\overline{H} \# K)_t = - H_t \circ \phi_H^t + K_t \circ \phi_H^t
= (K_t - H_t) \circ \phi_H^t . \tag 3.4 $$
Furthermore from the definitions of $\| \cdot \|$ and leng (see (1.3) and (1.4) respectively),
we have $\| H \| = \text{leng}\, (\phi_H)$.
In particular,
$$
\text{leng}(\phi_H^{-1}\phi_K) = \left\|\overline H \# K \right\|
= \| K - H \|. \tag 3.5
$$
The following simple lemma  will be useful later for the calculus
of the Hofer length function. The proof of this lemma immediately
follows from the definitions and is omitted.

\proclaim{Lemma 3.2} Let $H, \, K: [0,1] \times M \to \R$ be smooth.
Then we have \roster \item $\text{\rm leng}(\phi_H^{-1}\phi_K) = \text{\rm
leng}(\phi_K^{-1}\phi_H)$ \hskip0.5in or \hskip0.1in $\left\|\overline H\# K\right\| =
\left\|\overline K \# H\right\| = \| H - K \|$,
\smallskip

\item $ \text{\rm leng}(\phi_H \phi_K) \leq \text{\rm
leng}(\phi_H) + \text{\rm leng}(\phi_K)$ \hskip0.05in or \hskip0.1in $\|H \# K
\| \leq \|H\| + \|K\|,$
\smallskip

\item $\text{\rm leng}(\phi_H) = \text{\rm leng}(\phi_H^{-1})$
\hskip.9in or \hskip0.1in $\|H\| = \left\|\overline H \right\|.$
\endroster
\endproclaim

In relation to Floer homology and the spectral invariants, one
often needs to consider the periodic Hamiltonian functions $H$
satisfying $H(t+1,x) = H(t,x)$. For example, the spectral
invariants $\rho(\phi_H;a)$ of the Hamiltonian path $ \phi_H: t
\mapsto \phi_H^t$ are defined in [Oh4] first by reparameterizing
the path so that it becomes {\it boundary flat} (see Definition 3.3 below)
and so time-periodic in particular, by applying the
Floer homology theory to the Hamiltonian generating the
reparameterized Hamiltonian path, and then by proving the resulting
spectral invariants are independent of such reparameterization.
For this purpose, the senior author used the inequality
$$
\int_0^1 -\max(H-K)\,dt \leq \rho(\phi_H;a) - \rho(\phi_K;a) \leq
\int_0^1 - \min(H-K)\, dt
$$
in an essential way in [Oh4], [Oh5].

The following basic formula for the Hamiltonian generating a
reparameterized Hamiltonian path follows immediately from the definition.
It is used for the above purpose and again later in this paper.
For a given Hamiltonian function $H:\R \times M \to \R$,
not necessarily one-periodic, generating the Hamiltonian path
$\lambda = \phi_H$, the reparameterized path
$$
t \mapsto \phi_H^{\zeta(t)}
$$
is generated by the Hamiltonian function $H^\zeta$ defined by
$$
H^\zeta(t,x): = \zeta'(t)H(\zeta(t),x)
$$
for any smooth function $\zeta:\R \to \R$.
Here $\zeta'$ denotes the derivative of the function $\zeta$.
In relation to the reparameterization of Hamiltonian paths, the
following definition will be useful.

\definition{Definition 3.3} We call a path $\lambda:[0,1] \to
Symp(M,\omega)$ {\it boundary flat near $0$ (near $1$)} if $\lambda$
is constant near
$t=0$ ($t=1$), and we call the path {\it boundary flat}
if it is constant near $t=0$ and $t=1$.
\enddefinition

Of course this is the same as saying that any generating
Hamiltonian $H$ of $\lambda$ is constant near the end points. We would
like to point out that {\it the set of boundary flat Hamiltonians
is closed under the operations of the product $(H,K)\mapsto H\# K$
and taking the inverse $H \mapsto \overline H$} (and similarly for
paths that are flat near $t=0$ or $t=1$).

We will see in the $L^{(1,\infty)}$-Approximation Lemma (Appendix
2) that by choosing a suitable $\zeta$ so that $\zeta'\equiv 0$ near
$t=0, \, 1$ any Hamiltonian path can be approximated by a boundary
flat one in the Hamiltonian topology which we will introduce
later. We would like to emphasize that this approximation cannot be done in
the $L^\infty$-norm and that there is no such approximation procedure
in the $L^\infty$-topology. This would obstruct the smoothing
procedure of concatenated Hamiltonian paths or the extension of
the spectral invariants to the $C^0$-category (see [Oh7]), which is
the main reason why we adopt the $L^{(1,\infty)}$-norm, in addition to
its natural appearance in Floer theory.

Let $\lambda: [0,1] \to Symp(M,\omega)$ be a smooth path such that
$\lambda(t) \in Ham(M,\omega) \subset Symp(M,\omega)$.
We know that by definition of $Ham(M,\omega)$, for each given $s
\in [0,1]$ there exists a unique normalized Hamiltonian
$H^s=\{H^s_t\}_{0 \leq t \leq 1}$ such that $H^s \mapsto
\lambda(s)$. One very important property of a $C^\infty$-path (or
$C^1$ path in general) $\lambda: [0,1] \to Ham(M,\omega)$ is the following result by
Banyaga [Ba]
\proclaim{Proposition 3.4 [Proposition II.3.3, Ba]}
Let $\lambda: [0,1] \to Symp(M,\omega)$ be a smooth path such that
$\lambda(t) \in Ham(M,\omega) \subset Symp(M,\omega)$.
Define the vector field $\dot \lambda$ by
$$
\dot{\lambda}(s) := {\part \lambda \over \part s} \circ (\lambda(s))^{-1}
$$
and consider the closed one-form $\dot{\lambda}\rfloor \omega$.
Then this one-form is exact for all $s \in [0,1]$.
\endproclaim

In other words, {\it any smooth path in $Symp(M,\omega)$
whose image lies in $Ham(M,\omega)$ is
Hamiltonian} in the sense of Definition 1.2. Note that this
statement does not make sense if the path is not at least $C^1$ in $s$, i.e.,
when we consider a {\it continuous} path in $Homeo(M)$
whose image lies in $Ham(M,\omega)$.
As far as we know, it is not known whether one can always
approximate a continuous path $\lambda: [0,1] \to Ham(M,\omega)
\subset  Symp_0(M, \omega)\hookrightarrow Homeo(M)$ by a sequence
of smooth Hamiltonian paths. More precisely, it is not known in
general whether there is a
sequence of smooth Hamiltonian functions $H_j: [0,1] \times
M \to \R $ such that the Hamiltonian paths $t \mapsto
\phi^t_{H_j}$ uniformly converge to $\lambda$.

Not only for its definition but also for many results in the study
of the geometry of the Hamiltonian diffeomorphism group, a path
being Hamiltonian, not just lying in $Ham(M,\omega)$, is a crucial
ingredient. For that reason, it is reasonable to attempt to keep
track of the former property as one develops the {\it topological
Hamiltonian geometry}. Our definition of the Hamiltonian topology
in the present paper is the outcome of this attempt.

Obviously there is a one-one correspondence between the set of Hamiltonian
paths and that of generating (normalized) Hamiltonians in the smooth category.
However this correspondence gets murkier as the regularity of the
Hamiltonian gets worse, say when the regularity is less than
$C^{1,1}$. Because of this, we introduce the following terminology
for our later discussions.

\definition{Definition 3.5} We recall that $\PP^{ham}(Symp(M,\omega),id)$ denotes the
set of (smooth) Hamiltonian paths $\lambda$ defined on $[0,1]$
satisfying $\lambda(0) = id$ (see Definition 1.2 and Definition A.1).
Let $H$ be the (unique normalized) Hamiltonian
generating a given Hamiltonian path $\lambda$. We define two
maps
$$
\text{Tan}, \, \text{Dev} : \PP^{ham}(Symp(M,\omega),id) \to C_m^\infty([0,1]\times M, \R)
$$
by the formulas
$$
\aligned \text{Tan}(\lambda)(t,x)& := H(t,(\phi_H^t)(x)), \\
\text{Dev}(\lambda)(t,x) &: = H(t,x),
\endaligned
$$
and call them the {\it tangent map} and the {\it developing map}.
We call the image of the tangent map $\text{Tan}$ the {\it rolled
Hamiltonian} of $\lambda$ (or of $H$).
\enddefinition
The identity (3.2) implies the identity
$$
\text{Tan}(\lambda) = - \text{Dev}(\lambda^{-1}) \tag 3.6
$$
for a general (smooth) Hamiltonian path $\lambda$.

The tangent map corresponds to the map of the tangent vectors of the
path. Assigning the usual generating Hamiltonian $H$ to a
Hamiltonian path corresponds to the {\it developing map} in the Lie
group theory: one can `develop' any differentiable path on a Lie
group to a path in its Lie algebra using the tangent map and then by
right translation. (The senior author would like to take this
opportunity to thank A. Weinstein for making this remark almost 9
years ago right after he wrote his first papers [Oh1,2] on the
spectral invariants. Weinstein's remark answered the questions about
the group structure $(\#, -)$ on the space of Hamiltonians and much
helped the senior author's understanding of the group structure at
that time.)

We also consider the evaluation map
$$
ev_1: \PP^{ham}(Symp(M,\omega),id) \to Symp(M,\omega), \quad
ev_1(\lambda) = \lambda(1),
$$
and the obvious composition of maps
$$
\iota_{ham}:\PP^{ham}(Symp(M,\omega),id) \hookrightarrow
\PP(Symp(M,\omega),id) \to \PP (Homeo (M),id).
$$
We next state the following proposition. This proposition is a
reformulation of Theorem 6, Chapter 5 [HZ], in our general context,
which Hofer and Zehnder proved for compactly supported
Hamiltonian diffeomorphisms on $\R^{2n}$. In the presence of the
general energy-capacity inequality [LM], their proof can be easily adapted
to our general context. For readers' convenience, we give the
details of the proof here.

\proclaim{Proposition 3.6} Let $\lambda_i=\phi_{H_i} \in
\PP^{ham}(Symp(M,\omega),id)$ be a sequence of smooth Hamiltonian
paths and $\lambda=\phi_H$ be another smooth path such that
 \roster
\item $\left\|\overline H \# H_i\right\| \to 0$, and
\item $ev_1(\lambda_i) = \phi_{H_i}^1 \to \psi$ uniformly to a
map $\psi:M\to M$.
\endroster
Then we must have $\psi = \phi_H^1$.
\endproclaim
\demo{Proof} We first note that $\psi$ must be continuous since it
is a uniform limit of continuous maps $\phi_{H_i}^{1}$. Suppose
the contrary that $\psi \neq \phi^1_H$, i.e., $(\phi^1_H)^{-1}\psi
\neq id$. Then we can find a small closed ball $B$ such that
$$
B \cap \left( (\phi^1_H)^{-1}\psi \right)(B) = \emptyset.
$$
Since $B$ and hence $\left( (\phi^1_H)^{-1}\psi \right)(B)$
is compact and $\phi^1_{H_i}
\to \psi$ uniformly, we have
$$
B \cap \Big((\phi^1_H)^{-1}\phi_{H_i}^1\Big)(B)= \emptyset
$$
for all sufficiently large $i$. By definition of the Hofer
displacement energy $e$ (see [H1] for the definition), we have $e(B) \leq
\|(\phi^1_H)^{-1}\phi_{H_i}^1\|$. Now by the energy-capacity
inequality from [LM], we know $e(B) > 0$  and hence
$$
0 < e(B) \leq \left\|(\phi^1_H)^{-1}\phi_{H_i}^1\right\|
$$
for all sufficiently large $i$. On the other hand, we have
$$
\left\|(\phi^1_H)^{-1}\phi_{H_i}^1\right\| \leq \left\|\overline H \# H_i \right\| \to 0
$$
by hypothesis (1). The last two inequalities certainly
contradict each other.
That completes the proof. \qed \enddemo

What this proposition indicates for the practical purpose is that
simultaneously imposing both convergence
$$
\align & \|\overline H \# H_i \|  \to 0
\quad\text{and } \\
& \phi_{H_i}^1  \to \phi_{H}^1 \quad \text{in the $C^0$-topology}
\endalign
$$
is consistent in that they give rise to a nontrivial topology.

Remark that the evaluation map $ev_1$ is {\it not} continuous if
we equip $\PP^{ham}(Symp(M,\omega),id)$ with the Hofer topology
(Definition 1.3) and $Ham(M,\omega)$ with the $C^0$-topology (and
therefore Proposition 3.6 is not trivial). If it were, for every
sequence $H_i$ such that $\| H_i \| \to 0$, we would have
$\phi_{H_i}^1 \to id$. But, for any pair $(x,y)$ of points $x,\, y
\in M$, it is well-known that there is such a sequence with
$\phi_{H_i}^1 (x) = y$ for all $i$: This is because {\it the
transport energy of a point from one place to any other place is
always zero}, that is
$$
\inf_{H} \{\|H\| \mid \phi_H^1(x) = y\} = 0.
$$

\medskip

We will now define the (strong) Hamiltonian topology.
Its definition is directly motivated by the above Propositions 3.1 and 3.6
(see the remarks after the propositions).

\definition{Definition 3.7 [(Strong) Hamiltonian topology]}
\roster \item We define the {\it (strong) Hamiltonian topology} on
the set $\PP^{ham}(Symp(M,\omega),id)$ of Hamiltonian paths  by
the one generated by the collection of subsets
$$
\aligned & \UU(\phi_H,\e_1,\e_2) := \\
& \Big\{\phi_{H'} \in \PP^{ham}(Symp(M,\omega),id) \Big|
\|\overline H \# H' \| < \e_1, \, \overline d(\phi_H, \phi_{H'}) <
\e_2 \Big\}
\endaligned
\tag 3.7
$$
of $\PP^{ham}(Symp(M,\omega),id)$ for $\e_1, \, \e_2> 0 $ and
$\phi_H \in \PP^{ham}(Symp(M,\omega),id)$. We denote the resulting
topological space by $\PP_s^{ham}(Symp(M,\omega),id)$. \item We
define the {\it (strong) Hamiltonian topology} on $Ham(M,\omega)$ to be
the strongest topology such that the evaluation map
$$
ev_1: \PP_s^{ham}(Symp(M,\omega),id) \to Ham(M)
$$
is continuous. We denote the resulting topological space by $\HH
am(M,\omega)$.
\endroster
We will call continuous maps  with respect to the (strong)
Hamiltonian topology {\it (strongly) Hamiltonian continuous}.
\enddefinition

We refer readers to section 6 for the corresponding definition of
Hamiltonian topology  either for the non-compact case or the case
of manifolds with boundary.

We should now make several remarks concerning our choice of the
above definition of the Hamiltonian topology. The
combination of the Hofer topology and the $C^0$-topology in (3.7)
will be crucial to carry out all of the limiting process towards
the $C^0$-Hamiltonian world in this paper and in [Oh7]. Such a
phenomenon was first indicated by Eliashberg [El] and partly
demonstrated by Viterbo [V] and Hofer [H1,2].

We have the following interpretation of the Hamiltonian topology, which will be used later.

By definition, we have the natural continuous maps
$$
\aligned
\iota_{ham} \colon \PP_s^{ham}(Symp(M,\omega),id) & \to
\PP(Symp(M,\omega),id)
\hookrightarrow \PP(Homeo(M),id), \\
\text{\rm Dev} \colon \PP_s^{ham}(Symp(M,\omega),id) & \to
C^\infty_m([0,1]\times M, \R) \hookrightarrow
L_m^{(1,\infty)}([0,1]\times M, \R).
\endaligned
\tag 3.8
$$
We call the product map
$$
(\iota_{ham},\text{Dev}):\PP_s^{ham}(Symp(M,\omega),id) \to \PP(Symp(M,\omega),id)
\times C^\infty_m([0,1] \times M,\R)
$$
the {\it unfolding map}.
The Hamiltonian topology on
$\PP^{ham}(Symp(M,\omega),id)$ is nothing but the weakest topology for
which this unfolding map is continuous.

Here are several other comments.

\definition{Remark 3.8} \roster
\item The way how we define a topology on $Ham(M,\omega)$ starting
from one on the path space $\PP^{ham}(Symp(M,\omega),id)$ is
natural since the group $Ham(M,\omega)$ itself is defined that
way. We will repeatedly use this strategy in this paper.
\item Note that the Hamiltonian topology on $Ham (M,\omega)$
is nothing but the one induced by the evaluation map $ev_1$.
\item We also note that the collection of sets (3.7) is symmetric with respect to $H$ and $H'$,
i.e., $\phi_{H'} \in \UU (\phi_H, \e_1, \e_2 ) \iff \phi_H \in \UU (\phi_{H'}, \e_1, \e_2 )$.
\item It is easy to see that for fixed $\phi_H \in \PP^{ham}_s (Symp(M,\omega),id)$,
the open sets (3.7) form a neighborhood basis of the Hamiltonian topology at $\phi_H$.
\item Because of the simple identity
$$
(\overline H \# H')(t,x) = (H' - H)(t, \phi_H^t(x))
$$
one can write the length in either of the following two ways:
$$
\text{leng}(\phi_H^{-1} \phi_{H'}) = \| \overline H \# H' \| =
\|H' - H\|
$$
if $H$ and $H'$ are smooth (or more generally $C^{1,1}$). In this
paper, we will mostly use the first one that manifests the group
structure better. The proof is straightforward to check and
omitted.
\item
Note that the above identity does not make sense in general even
for $C^1$-functions because their Hamiltonian vector field would
be only $C^0$ and so their flow $\phi_H^t$ may not exist.
Understanding what is going on in such a case touches the heart of
the $C^0$-Hamiltonian geometry and dynamics. We will pursue the
dynamical issue in [Oh7] and focus on the geometry in this paper.
\endroster
\enddefinition

It turns out that $\PP_s^{ham}(Symp(M,\omega),id)$ is metrizable.
We now define the following natural metric on $\PP^{ham}(Symp(M,\omega),id)$
which combines the Hofer metric and the $C^0$-metric appropriately.

\definition{Definition 3.9}
We define a metric on $\PP^{ham}(Symp(M,\omega),id)$ by
$$ d_{ham}(\phi_H,\phi_{H'}) =
\|\overline H \# H' \| + \overline d(\phi_H, \phi_{H'}). $$
\enddefinition

\proclaim{Proposition 3.10}
The Hamiltonian topology on $\PP^{ham}(Symp(M,\omega),id)$ is equivalent to
the metric topology induced by $d_{ham}$.
\endproclaim

\demo{Proof} This is an exercise in using the definitions.
Let $\UU$ be open in the Hamiltonian topology, and let $\phi_H \in \UU$.
By Remark 3.8(4), there are $\e_1, \e_2 > 0$ such that $\UU (\phi_H, \e_1, \e_2) \subset \UU$.
Define $\e = \min (\e_1, \e_2)$.
Let
$$
\UU_\e (\phi_H) =\{ \phi_{H'} \in \PP_s^{ham}(Symp(M,\omega),id)
\mid d_{ham} (\phi_H, \phi_{H'}) < \e \}
$$
be the metric ball of radius $\e$ centered at $\phi_H$.
By our choice for $\e$ and by Definitions 3.7(1) and 3.9,
we have $\UU_\e (\phi_H) \subset \UU (\phi_H, \e_1, \e_2) \subset \UU$.
This holds for any $\phi_H \in \UU$, so $\UU$ is open in the metric topology.

Conversely, suppose $\VV$ is open in the metric topology, and $\phi_H \in \VV$.
Then $\UU_\e (\phi_H) \subset \VV$ for some $\e > 0$, and
$\UU (\phi_H, \frac{\e}{2}, \frac{\e}{2}) \subset \UU_\e (\phi_H) \subset \VV$.
So $\VV$ is open in the metric topology.
\qed\enddemo

\proclaim{Proposition 3.11} The left translations of the
group $\PP^{ham}_s(Symp(M,\omega),id)$ are continuous, i.e.,
for each $\lambda \in \PP^{ham}_s(Symp(M,\omega),id)$, the bijection
$$
L_\lambda:\PP_s^{ham}(Symp(M,\omega),id) \to \PP_s^{ham}(Symp(M,\omega),id),
\quad L_\lambda(\mu) = \lambda\mu,
$$
is continuous, and therefore a homeomorphism, with respect to the Hamiltonian topology
on $\PP^{ham}_s(Symp(M,\omega),id)$. In particular, the sets of the form
$$
\phi_H \left( \UU (id, \e_1, \e_2 ) \right), \quad \e_1, \, \e_2 > 0
\tag 3.9
$$
form a neighborhood basis at $\phi_H$ in $\PP^{ham}_s (Symp(M,\omega),id)$.
\endproclaim
\demo{Proof} Let $\lambda = \phi_H$. We have to show that
$L_\lambda^{-1}(\UU (\phi_K, \e_1, \e_2))$ is open for any choice
of $\mu = \phi_K$ and $\e_1, \, \e_2 > 0$. Let $\phi_L \in
L_\lambda^{-1}(\UU (\phi_K, \e_1, \e_2))$, i.e.,
$$
\phi_H \phi_L \in \UU (\phi_K, \e_1, \e_2). \tag 3.10
$$
We need to find some $\e_1', \, \e_2' > 0$ such that
$$
\UU(\phi_L,\e_1',\e_2') \subset L_\lambda^{-1}(\UU (\phi_K, \e_1,
\e_2)),
$$
or equivalently, such that
$$
L_\lambda(\UU(\phi_L,\e_1',\e_2')) =
\phi_H(\UU(\phi_L,\e_1',\e_2')) \subset \UU (\phi_K, \e_1, \e_2).
\tag 3.11
$$

For the part of $\overline d$, we define
$$
\bar \e_2 = \e_2 - \overline d(\phi_H\phi_L,\phi_K) > 0 \tag 3.12
$$
by (3.10).
By compactness of $M$, the smooth map $[0,1] \times M \to M, (t,x) \mapsto \phi_H^t (x)$
is in particular uniformly continuous with respect to the standard metric on $[0,1]$
and the metric $d$ on $M$.
Therefore there exists $0<\e_2'<\bar\e_2$ such that
$$
d(x,y) < \e_2' \quad \Longrightarrow \quad d (\phi_H^t(x), \phi_H^t(y)) < \bar
\e_2
$$
for all $x,y \in M$ and all $t \in [0,1]$. Hence if $\overline d(\phi_L, \phi_{L'}) <
\e_2'$, then
$$
\aligned \overline d (\phi_H \phi_L, \phi_H \phi_{L'}) & =
\max\{d_{C^0}(\phi_H \phi_L, \phi_H \phi_{L'}),
d_{C^0} (\phi_L^{-1} \phi_H^{-1}, \phi_{L'}^{-1} \phi_H^{-1})\} \\
& = \max \left\{ \max_{(t,x)} d \left(\phi_H^t \phi_L^t (x), \phi_H^t \phi_{L'}^t (x) \right),
d_{C^0} ( \phi_L^{-1}, \phi_{L'}^{-1} ) \right\} \\
& < \max \{\bar \e_2, \e_2' \} = \bar \e_2.
\endaligned
$$
We now estimate
$$
\overline d(\phi_H\phi_{L'},\phi_K) \leq \overline d (\phi_H
\phi_{L'}, \phi_H \phi_L) + \overline d(\phi_H \phi_L, \phi_K) < \bar \e_2 +
\overline d(\phi_H \phi_L, \phi_K) = \e_2 \tag 3.13
$$
by (3.12), as long as $\overline d(\phi_L, \phi_{L'}) < \e_2'$.

On the other hand for the part of $\|\cdot\|$,
choose $\e_1' = \e_1 - \|H \# L - K\|$, which again is positive by (3.10).
It is immediate to check from the definitions that $\|H \# L' - H\# L\| = \|L' - L\|$.
Then whenever $L'$ satisfies $\|L'-L\| < \e_1'$,
we have by the triangle inequality
$$
\|H\# L' - K\| \leq \|H \# L' - H\# L\| + \|H \#L -K\| = \|L' - L\| + \|H \# L - K\| < \e_1.
$$

That completes the proof of the first statement.
Since the inverse of $L_\lambda$ is the left translation $L_{\lambda^{-1}}$,
left translations are in fact homeomorphisms.
The last statement is obvious from this and Remark 3.8(4).
This finishes the proof. \qed\enddemo

As we will see below, $\PP_s^{ham}(Symp(M,\omega),id)$ in fact forms a topological group.
This will follow as a corollary to the fact that its completion $\overline{\PP_s^{ham}(Symp(M,\omega),id)}$
considered below forms a topological group as well.
But we prefer to give an elementary proof of Proposition 3.11 and the following corollaries
using only the definitions,
and then to complete the discussion of $\PP_s^{ham}(Symp(M,\omega),id)$ and $\HH am(M,\omega)$,
before dealing with the more complicated arguments involved when considering said completion.

Proposition 3.11 immediately gives rise to the following corollaries.

\proclaim{Corollary 3.12} The evaluation map $ev_1:\PP^{ham}_s(Symp(M,\omega),id) \to \HH am(M,\omega)$
is an open map with respect to the Hamiltonian topology on $Ham(M,\omega)$.
In particular, the following hold:
\roster
\item For fixed $\phi \in Ham(M,\omega)$ and $H \mapsto \phi$,
the sets of the form
$$
ev_1 \Big(\UU(\phi_H, \e_1, \e_2) \Big), \quad \e_1, \, \e_2 > 0
$$
form a neighborhood basis at $\phi$ in the Hamiltonian topology.
\item For fixed $\phi \in Ham(M,\omega)$ and $H \mapsto \phi$,
the sets of the form
$$ \phi \Big( ev_1 \Big( \UU (id, \e_1, \e_2) \Big) \Big) =
ev_1 \Big( \phi_H \Big( \UU (id, \e_1, \e_2) \Big) \Big), \quad \e_1, \, \e_2 > 0
$$
also form a neighborhood basis at $\phi$ in the Hamiltonian topology.
\endroster
\endproclaim

\demo{Proof} Let $\UU \subset \PP^{ham}_s(Symp(M,\omega),id)$ be open in the Hamiltonian topology.
We have to show that $ev_1 (\UU) \subset \HH am(M,\omega)$
is open with respect to the Hamiltonian topology on $Ham (M,\omega)$.
But by definition of the Hamiltonian topology, $ev_1(\UU)$ is open if and only if
$$ ev_1^{-1} \left(ev_1 (\UU)\right) =
\bigcup_\lambda \{ \lambda(\UU)\mid \lambda \in
\PP^{ham}_s(Symp(M,\omega),id), \lambda(0) = \lambda (1) = id \} $$
is open. But the latter is the union of open sets by Proposition
3.11 and hence itself open. That proves the first part.

Openness and continuity of $ev_1$ with respect to the Hamiltonian topology
together with Remark 3.8(4) now implies (1).

For (2), note that since $Ham(M,\omega)$ is a group
it also acts on itself via left translations.
The left translations of $\PP^{ham}_s(Symp(M,\omega),id)$
and $Ham(M,\omega)$ commute with $ev_1$ in the sense that
if $\phi \in Ham(M,\omega)$ and $H \mapsto \phi$ is any Hamiltonian, then
$ev_1 (\phi_H \phi_{H'}) = \phi (ev_1 (\phi_{H'}))$ for any
$\phi_{H'} \in \PP^{ham}_s (Symp(M,\omega),id)$.
In other words, $ev_1$ is a (continuous) group homomorphism.
This together with openness and continuity of $ev_1$ and
the last statement of Proposition 3.11 implies (2). \qed
\enddemo

The following is one indication of good properties of the Hamiltonian
topology.

\proclaim{Theorem 3.13} $\HH am(M,\omega)$ is path-connected and locally path-connected.
\endproclaim

\demo{Proof} We first prove that $\HH am(M,\omega)$ is locally path-connected at the identity.
Consider the following open neighborhood of the identity element in $\HH am(M,\omega)$
$$
\UU = ev_1 \Big( \UU (id, \e_1, \e_2) \Big)
$$
for any $\e_1, \, \e_2 > 0$.
Note that by Corollary 3.12 these sets form a neighborhood basis at the identity.
So it suffices to prove that $\UU$ is path-connected.

Let $\phi_0 \in \UU$.
We will prove that $\phi_0$ can be connected by a continuous
path to the identity inside $\UU$.
Since $\phi_0 \in \UU$ there exists $H \mapsto \phi_0$ such that
$$
\|H\| < \e_1, \quad
\overline d (\phi_H, id) = \sup_{t\in [0,1]} \overline d(\phi_H^t, id) < \e_2 .
$$
Let $H^s$ be the Hamiltonian generating $t \mapsto \phi_{H^s}^t = \phi_H^{st}$
defined by $H^s(t,x) = s H(st,x)$. We have
$$
\overline d (\phi_{H^s}, id) = \sup_{t \in [0,1]} \overline d(\phi_{H^s}^t,id) =
\sup_{t \in [0,s]} \overline d(\phi_H^t,id) \leq \sup_{t\in [0,1]} \overline d(\phi_H^t, id) < \e_2 .
$$
Also note that by substituting $\tau = st$ we get $\| H^s \| \leq \| H \|$.
Combining the two, we derive that $\phi_{H^s} \in \UU (id, \e_1, \e_2)$
and hence $\phi_H^s = \phi_{H^s}^1 \in \UU$ for all $s \in [0,1]$.
Hence the path $\lambda = \phi_H: t \mapsto \phi_H^t$
has its image contained in $\UU$, and connects the identity and $\phi_0$.
Continuity follows from Corollary A.3.
So $\UU$ is path-connected.

Now let $\phi \in Ham(M,\omega)$. By Corollary 3.12, the sets
$\phi \, \UU$, where $\UU$ as above, form a neighborhood basis at
$\phi$. That they are path-connected follows from their definition
and path-connectedness of $\UU$. This proves local
path-connectedness of $\HH am(M,\omega)$. Path-connectedness of
$\HH am(M,\omega)$ follows from its definition (see the remark
after Definition A.1) and Corollary A.3. That proves the theorem.
\qed\enddemo

\medskip

One crucial point of imposing the $C^0$-requirement in the Hamiltonian
topology compared to the Hofer topology is that it enables us to
extend the evaluation map $ev_1: \PP_s^{ham}(Symp(M,\omega),id) \to
Ham(M,\omega)$ to the completion of $\PP_s^{ham}(Symp(M,\omega),id)$
with respect to the corresponding metric topology.
Recall that the evaluation map is {\it not} continuous if one equips
$\PP^{ham}(Symp(M,\omega),id)$ with the Hofer topology and
$Ham(M,\omega)$ with the $C^0$-topology
(see the remark after Proposition 3.6).
It is also an interesting problem to understand the
completion of $Ham(M,\omega)$ with respect to the Hofer topology,
but this is much harder to study, partly because a general element
in the completion would not be a continuous map.

\medskip

We now define the notion of {\it topological Hamiltonian path}, {\it topological Hamiltonian function},
and {\it Hamiltonian homeomorphism}.
Let $(\phi_i, \lambda_i, H_i)$ be a sequence of triples,
where $\phi_i \in Ham (M,\omega)$ are Hamiltonian diffeomorphisms,
and $H_i \in C^\infty_m ([0,1] \times M, \R)$ are normalized Hamiltonian functions,
such that $H_i$ generates the Hamiltonian path $\lambda_i = \phi_{H_i} : t \mapsto \phi_{H_i}^t$,
and $\phi_i = \phi_{H_i}^1 = \lambda_i (1)$.
Suppose the sequence is Cauchy in the Hamiltonian topology,
$$
\align
\overline{d} \left( \phi_{H_i}, \phi_{H_j} \right) & \to 0, \quad \text{as } i, \, j \to \infty, \text{ and} \\
\| H_i - H_j \| & \to 0 , \quad \text{as } i, \, j \to \infty .
\endalign
$$
In particular, $H_i$ converges to a (normalized) $L^{(1,\infty)}$-function $H \in L^{(1,\infty)}_m ([0,1] \times M, \R)$,
$\lambda_i$ converges to a continuous path $\lambda \in \PP (Homeo (M), id)$, and
$\lambda (1) = \lim_{C^0} \phi_i =: h \in Homeo (M)$.
We call the continuous path $\lambda$ a {\it topological Hamiltonian path},
the function $H$ a {\it topological Hamiltonian function}, and the map $h$ a {\it Hamiltonian homeomorphism}.

More precisely, recall the unfolding map
$$
\align
(\iota_{ham},\text{Dev}) & :\PP_s^{ham}(Symp(M,\omega),id) \to
\PP(Symp(M,\omega),id) \times C_m^\infty([0,1] \times M,\R)\\
& \quad \quad \quad \hookrightarrow \PP(Homeo(M),id) \times
L_m^{(1,\infty)}([0,1] \times M,\R)
\endalign
$$
which was defined by $\lambda = \phi_H \mapsto (\lambda,H)$. We
denote by $\QQ$ the image of $(\iota_{ham},\text{Dev})$ equipped
with the subspace topology. More precisely, the topology on $\QQ$
is induced by the product metric given by the $C^0$-metric
$\overline d$ on $\PP(Homeo(M),id)$ and the
$L^{(1,\infty)}$-metric on $L_m^{(1,\infty)}([0,1]\times M, \R)$.
We will refer to this topology on $\QQ$ also as the Hamiltonian
topology. This will be further explained in Remark 3.17(2) below.

Note that Definition 3.9
implies that both $\iota_{ham}$ and $\text{Dev}$ are Lipschitz continuous
(with $L \leq 1$) with respect to $d_{ham}$ on $\PP_s^{ham}(Symp(M,\omega),id)$,
and the $C^0$-metric $\overline d$ on $\PP(Homeo(M),id)$ and the $L^{(1,\infty)}$-metric on
$L_m^{(1,\infty)}([0,1]\times M, \R)$ respectively.
These maps induce natural (Lipschitz continuous) projections from $\QQ$ onto the first and second factor,
denoted by
$$
\aligned
\iota_{ham}^Q &: \QQ \to \PP(Symp(M,\omega),id) \hookrightarrow \PP(Homeo(M),id), \\
Dev^Q & : \QQ \to C_m^\infty([0,1] \times M,\R)
\hookrightarrow L_m^{(1,\infty)}([0,1] \times M,\R).
\endaligned
\tag 3.14
$$
The map $ev_1$ is also seen to be Lipschitz continuous (also with $L \le 1$)
with respect to $d_{ham}$ on $\PP_s^{ham}(Symp(M,\omega),id)$ and
the $C^0$-topology on $Ham(M,\omega) \subset Homeo(M)$,
and hence induces the natural (Lipschitz continuous) map
$$
ev_1^Q : \QQ \to Ham(M,\omega) \subset Homeo(M), \quad (\lambda,H) \mapsto \lambda(1).
$$

We denote by $\overline{\QQ}$ the closure of $\QQ$ in
$\PP(Homeo(M),id) \times L^{(1,\infty)}_m([0,1] \times M,\R)$ with
respect to the product metric, and call any element thereof a {\it
strong Hamiltonian path}. By Lipschitz continuity of the above
maps, all three maps naturally extend to continuous maps defined
on $\overline{\QQ}$.

\definition{Definition 3.14 [(Strong) Hamiltonian homeomorphisms]}
We denote by
$$
\overline{ev}_1^Q : \overline{\QQ} \to Homeo(M), \quad
(\lambda,H) \mapsto \lambda(1) \tag 3.15
$$
the natural continuous extension of the evaluation map $ev_1^Q$.
We denote by
$$
Hameo(M,\omega) \subset Homeo(M)
$$
the image of $\overline{\QQ}$
under the map $\overline{ev}_1^Q$ and call any element thereof a {\it (strong) Hamiltonian homeomorphism}.
I.e., $h \in Hameo(M,\omega)$ if and only if there exists a Cauchy sequence
$(\phi_{H_i},H_i)$ in $\QQ$ in the Hamiltonian topology with
$h = \lim_{C^0}\phi_{H_i}^1$.
We equip $Hameo(M,\omega)$ with the subspace topology induced from $Homeo(M)$,
i.e., with the $C^0$-topology.
We define the {\it (strong) Hamiltonian topology} on the set $Hameo(M,\omega)$
to be the strongest topology such that the map $\overline{ev}_1^Q$ is continuous.
We denote by $\HH ameo(M,\omega)$ the resulting topological space.
By definition the map
$$
\overline{ev}_1^Q: \overline{\QQ} \to \HH ameo(M,\omega) \tag 3.16
$$
is surjective, continuous, and the following diagram commutes
$$
\matrix \QQ & \longrightarrow & \HH am(M,\omega) \\
\downarrow & \quad & \downarrow \\
\overline{\QQ} & \longrightarrow & \HH ameo(M,\omega),
\endmatrix \tag 3.17
$$
where the vertical maps are the natural inclusions,
and the horizontal maps are the maps induced by the evaluation map.
\enddefinition

\definition{Definition 3.15 [Topological Hamiltonian path]}
We denote by
$$
\overline{\iota}_{ham}^Q: \overline{\QQ} \to \PP(Homeo(M),id),
\quad (\lambda,H) \mapsto \lambda
$$
the natural continuous extension of the map $\iota_{ham}^Q$. By
the definition of $Sympeo(M,\omega)$ it follows that the map is
factorized into
$$
\overline{\iota}^Q_{ham}: \overline{\QQ} \to \PP(Sympeo(M,\omega),id) \hookrightarrow
\PP(Homeo(M),id).
$$
We denote by
$$
\PP^{ham}(Sympeo(M,\omega),id) \subset \PP(Sympeo(M,\omega),id)
$$
the image of the map $\overline{\iota}_{ham}^Q$ equipped
with the subspace topology, i.e., the $C^0$-topology.
We call any element $\lambda \in \PP^{ham}(Sympeo(M,\omega),id)$
a {\it topological Hamiltonian path}.
\enddefinition
\medskip

More specifically, a continuous path $\lambda \in
\PP(Homeo(M),id)$ is a topological Hamiltonian path if and only if
there exists a Cauchy sequence $(\phi_{H_i},H_i) \in \QQ$ in the
Hamiltonian topology such that $\lim_{C^0}\phi_{H_i} = \lambda$.

Now we ask the following uniqueness question on the
`$L^{(1,\infty)}$-Hamiltonian' concerning the one-oneness of the
map $\overline{\iota}_{ham}^Q$.

\definition{Question 3.16}
Consider the Cauchy sequences $(\phi_{H_i},H_i)$ and $(\phi_{H_i'},H_i')$
in the Hamiltonian topology such that
$(\phi_{H_i}^t)^{-1} (\phi_{H_i'}^t) \to id$
as $i \to \infty$ uniformly over $[0,1] \times M$.
Does this imply $\left\|\overline H_i \# H_i'\right\| \to 0$ as $i \to \infty$?
\enddefinition

The $C^0$-(or $L^\infty$-)version of this question
has been answered affirmatively by Viterbo [V2],
and then subsequently in the $L^{(1,\infty)}$-case
by the senior author [Oh7] during the preparation of the current revision of the paper.
We refer readers to [Oh7] for the generalization of this uniqueness
result in the Lagrangian context and for several other consequences
of this uniqueness result.

Here are several remarks.

\definition{Remark 3.17}
\roster \item Similarly, we can define the continuous extension
$\overline{Dev^Q}$ of $Dev^Q$. The image of this map is by
definition the set of topological Hamiltonian functions. These
will be studied in a sequel [Oh7]. \item Of course, as topological
spaces $\PP_s^{ham}(Symp(M,\omega),id) \cong \QQ$ via the
unfolding map. But it is often more convenient to consider the
completion of $\QQ$ in $\PP(Homeo(M),id) \times
L^{(1,\infty)}_m([0,1] \times M,\R)$ rather than the abstract
completion $\overline{\PP_s^{ham}(Symp(M,\omega),id)}$ of
$\PP_s^{ham}(Symp(M,\omega),id)$, and then dealing with
equivalence classes of Cauchy sequences representing elements in
$\overline{\PP_s^{ham}(Symp(M,\omega),id)}$. As topological spaces
$\overline{\PP_s^{ham}(Symp(M,\omega),id)} \cong \overline{\QQ}$
via the natural extension of the unfolding map. All statements
about $\QQ$ and $\overline{\QQ}$ can be translated to
$\PP_s^{ham}(Symp(M,\omega),id)$ and
$\overline{\PP_s^{ham}(Symp(M,\omega),id)}$ by composing all maps
with the unfolding map or its inverse, and vice versa. \item The
way how we define $Hameo(M,\omega)$ starting from the completion
of the path space $\PP_s^{ham}(Symp(M,\omega),id)$ is natural
since $Ham(M,\omega)$ itself is defined in a similar way (recall
Remark 3.8(1)).
\endroster
\enddefinition

Next recall $\text{Dev}(\phi_H)(t,x) = H(t,x)$ and $\text{Tan}(\phi_H)(t,x) = H(t, (\phi_H^t)(x))$.
For convenience, we will often write $H \circ \phi_H$ to denote
$(H \circ \phi_H) (t,x) = H (t, \phi_H^t (x)) = \text{Tan}(\phi_H)(t,x)$.
Note that from the definitions we immediately get the useful identity
$$
\text{leng} \Big( \phi_H (\phi_{H'})^{-1} \Big) = \| H \#
\overline{H'} \| = \| \text{Tan} (\phi_H) - \text{Tan} (\phi_{H'})
\|. \tag 3.18
$$

Continuity of the maps $\text{Dev}$ and $\text{Dev}^Q$ is obvious from their definition,
but not so that of $\text{Tan}$ and $\text{Tan}^Q$.
In this regard, we state the following lemma

\proclaim{Lemma 3.18} The map
$$
\text{Tan} : \PP_s^{ham}(Symp(M,\omega),id) \to C_m^\infty([0,1]\times M,\R)
$$
is continuous with respect to the
Hamiltonian topology on $\PP_s^{ham}(Symp(M,\omega),id)$ and the
$L^{(1,\infty)}$-topology on $C_m^\infty([0,1]\times M,\R)$.
The same holds for the map
$$
\text{Tan}^Q : \QQ \to C_m^\infty([0,1]\times M,\R), \quad (\lambda,H) \mapsto H \circ \lambda.
$$
\endproclaim
\demo{Proof}
Let $\lambda = \phi_{H}$ be given.
Consider another Hamiltonian path $\lambda' = \phi_{H'}$. We have
$$
\align \| \text{Tan}(\phi_{H'}) - \text{Tan}(\phi_H) \| & =
\| H' \circ \phi_{H'} - H \circ \phi_{H} \| \\
& \leq \| H' \circ \phi_{H'} - H \circ \phi_{H'} \| +
\|H \circ \phi_{H'} - H\circ \phi_H \| \\
& \leq \| H' - H \| + 2 \cdot L \cdot d_{C^0} (\phi_{H'}, \phi_H),
\tag 3.19
\endalign
$$
where $L$ is a Lipschitz constant that depends only on the smooth function $H$.
It follows from this inequality that $\text{Tan}$ is continuous
at every $\lambda \in \PP_s^{ham}(Symp(M,\omega),id)$ and hence the proof.
The proof for $\text{Tan}^Q$ is of course the same.
\qed \enddemo

Since the constant $L$ in (3.19) depends on the
Hamiltonian function $H$, the map $\text{Tan}$ is unlikely
to be {\it uniformly continuous}.
The constant $L$ cannot be controlled in the Hamiltonian topology,
e.g., when we consider a Cauchy sequence $(\phi_{H_i},H_i)$
representing a strong Hamiltonian path. This was the source
of many erroneous statements and proofs in the previous senior
author's own versions of the current paper, many of which are
corrected by the junior author in the current version.
The crucial lemma to deal with this difficulty will be
the Reparameterization Lemma 3.21 below.

Very often in the study of the geometry of Hamiltonian
diffeomorphisms, one needs to reparameterize a given Hamiltonian
path in a way that the reparameterization is close enough to the
given parameterization, e.g., in the smoothing process of the
concatenation of two paths. We now provide the correct topology
describing the closeness of such parameterizations.

\definition{Definition 3.19} We call the norm
$$
\|f\|_{ham}: = \|f\|_{C^0} + \|f'\|_{L^1}
$$
of a (smooth) function $f : [0,1] \to \R$ the {\it hamiltonian norm}
of the function $f$.
Here $f'$ denotes the derivative of the function $f$.
We say that two smooth functions $\zeta_1, \,
\zeta_2: [0,1] \to [0,1]$ are {\it hamiltonian-close} to each other if
the norm
$$
\align
\|\zeta_1 - \zeta_2\|_{ham} & : = \|\zeta_1 - \zeta_2\|_{C^0}
+ \|\zeta_1' - \zeta_2'\|_{L^1}\\
& = \max_{t \in [0,1]} |\zeta_1(t) - \zeta_2(t)| + \int_0^1
|\zeta_1'(t) - \zeta_2'(t)|\, dt
\endalign
$$
is small.
\enddefinition

Recall that for a given Hamiltonian function $H$ generating the Hamiltonian path $\phi_H$,
the reparameterized path $t \mapsto \phi_H^{\zeta(t)}$
is generated by the Hamiltonian function $H^\zeta$ defined by
$H^\zeta (t,x) = \zeta'(t) H (\zeta(t),x)$,
where $\zeta'$ again denotes the derivative
of the reparameterization function $\zeta : [0,1] \to [0,1]$.

\proclaim{Lemma 3.20} Let $H : [0,1] \times M \to \R$ be a normalized smooth Hamiltonian function,
and let $\zeta_1$, $\zeta_2 : [0,1] \to [0,1]$ be two smooth reparameterization functions.
Then
$$ \| H^{\zeta_1} - H^{\zeta_2} \| \le C \| \zeta_1 - \zeta_2 \|_{ham}, \tag 3.20 $$
where $C \le 2 \max (\| H \|_{C^0}, L)$ is a constant that depends only on the $C^0$-norm
$$ \| H \|_{C^0} = \max_{(t,x)} | H (t,x) | < \infty $$
of $H$ and a Lipschitz constant (in the time variable) $L$ for $H$.
\endproclaim

We refer to Appendix 2 for the proof of Lemma 3.20.
But note that Lemma 3.20 does {\it not} hold if we replace the hamiltonian norm
by the $C^0$-norm of $\zeta_1 - \zeta_2$ in (3.20).

We now state the following useful lemma

\proclaim{Lemma 3.21 [Reparameterization Lemma]} Suppose $H_i :
[0,1] \times M \to \R$ is a Cauchy sequence of smooth functions in
the $L^{\left(1,\infty\right)}$-topology, i.e.,
$$ \| H_i - H_j \| \to 0 \quad \text{as} \quad i,j \to \infty, $$
$\zeta_1, \, \zeta_2: [0,1] \to [0,1]$ are smooth
reparameterization functions on $[0,1]$, and $\lambda$, $\mu \in
\PP (Homeo(M),id)$ are continuous paths. Let $\e > 0$ be given.
\roster \item Then there exist $\delta = \delta(H_i) > 0$ and $i_0
= i_0 (H_i) > 0$ such that
$$ \| H_i^{\zeta_1} - H_i^{\zeta_2} \| < \e $$
for all $i \geq i_0$,
if $\zeta_1, \, \zeta_2$ satisfy
$$ \| \zeta_1 - \zeta_2 \|_{ham} < \delta . $$
\item
There exist $\delta' = \delta'(H_i) > 0$ and $i'_0 = i'_0 (H_i) > 0$ such that
$$ \| H_i \circ \lambda - H_i \circ \mu  \| < \e $$
for all $i \geq i'_0$,
if $\lambda, \, \mu$ satisfy
$$ d_{C^0} ( \lambda, \mu) < \delta' .
$$
\endroster
\endproclaim

\demo{Proof} (1) We can find $i_0$ sufficiently large such that
$$ \| H_i - H_{i_0} \| < \frac{\e}{3} \quad \text{for all} \quad i \ge i_0 . $$
Choose $0 < \delta < \frac{\e}{3 C}$, where $C$ is as in Lemma 3.20 with $H$ replaced by $H_{i_0}$. Then
$$ \| H_{i_0}^{\zeta_1} - H_{i_0}^{\zeta_2} \| < \frac{\e}{3} \quad \text{when} \quad
\|\zeta_1 - \zeta_2\|_{ham} < \delta.
$$
Therefore
$$
\align \| H_i^{\zeta_1} - H_i^{\zeta_2} \| & \leq \| H_i^{\zeta_1}
- H_{i_0}^{\zeta_1} \| + \| H_{i_0}^{\zeta_1} - H_{i_0}^{\zeta_2}
\| +
\| H_{i_0}^{\zeta_2} - H_i^{\zeta_2} \| \\
& = \| H_i - H_{i_0} \| + \| H_{i_0}^{\zeta_1} - H_{i_0}^{\zeta_2}
\| +
\| H_{i_0} - H_i \| \\
& < \frac{\e}{3} + \frac{\e}{3} + \frac{\e}{3} = \e ,
\endalign
$$
when $\|\zeta_1 - \zeta_2\|_{ham} < \delta$, $i \geq i_0$. That proves (1).

For (2), again choose $i'_0 = i_0$ sufficiently large such that
$$ \| H_i - H_{i_0} \| < \frac{\e}{3} \quad \text{for all} \quad i \ge i_0 . $$
By uniform continuity of $H_{i_0}$ there exists $\delta' > 0$ such that
$$ \| H_{i_0} \circ \lambda - H_{i_0} \circ \mu \|_\infty < \frac{\e}{6} $$
when $d_{C^0} (\lambda, \mu) < \delta$.
This implies
$$ \| H_{i_0} \circ \lambda - H_{i_0} \circ \mu \| < \frac{\e}{3} $$
when $d_{C^0} (\lambda, \mu) < \delta$.
Now apply the triangle inequality as above. \qed\enddemo

Note that $H_i$ converges to an $L^{(1,\infty)}$-function $H$,
but that we cannot replace $H_{i_0}$ by $H$ in the above proof
since $H$ is not even continuous in general.

\proclaim{Proposition 3.22} There exist continuous maps
$\overline{\text{Tan}^Q}$ and $\overline{\text{Dev}^Q}$, which we
again call the {\it tangent map} and the {\it developing map}
respectively
$$
\overline{\text{Tan}^Q}, \, \overline{\text{Dev}^Q}: \overline{\QQ}
\to L_m^{(1,\infty)}([0,1]\times M), \tag 3.21
$$
such that the following diagram commutes
$$
\matrix \QQ & \longrightarrow \quad & C_m^\infty([0,1] \times M, \R) \\
\downarrow & \quad & \downarrow \\
\overline{\QQ} & \longrightarrow \quad & L_m^{(1,\infty)}([0,1]\times M, \R),
\endmatrix \tag 3.22
$$
where the vertical maps are the natural inclusions,
and the horizontal maps are the tangent and developing maps.
\endproclaim

\demo{Proof} $\overline{\text{Dev}^Q}$ was already considered
above. For $\overline{\text{Tan}^Q}$, recall that
$$
\align \| \text{Tan}(\phi_{H_i}) - \text{Tan}(\phi_{H_j}) \| & =
\| H_i \circ \phi_{H_i} - H_j \circ \phi_{H_j} \| \\
& \leq \| H_i \circ \phi_{H_i} - H_j \circ \phi_{H_i} \|
+ \| H_j \circ \phi_{H_i} - H_j \circ \phi_{H_j} \| \\
& = \| H_i - H_j \| + \| H_j \circ \phi_{H_i} - H_j \circ \phi_{H_j} \|.
\endalign
$$

Now if $(\phi_{H_i}, H_i)$ is a Cauchy sequence in the Hamiltonian topology,
then the first term converges to zero by definition,
and the second term converges by Lemma 3.21(2).
So $\text{Tan}(\phi_{H_i})$ converges to an element in
$L_m^{(1,\infty)}([0,1]\times M, \R)$.

If $(\lambda,H) \in \overline{\QQ}$, there exists such a Cauchy sequence
$(\phi_{H_i}, H_i)$ converging to $(\lambda,H)$ in the Hamiltonian topology.
By definition,
$$
\overline{\text{Tan}^Q} (\lambda,H) = \lim_{i \to \infty}
\text{Tan} (\phi_{H_i}) = H\circ \lambda.
$$
Here the composition $H\circ \lambda$ is already well-defined as
an $L^{(1,\infty)}$-function.

Now suppose $(\lambda,H) \in \overline{\QQ}$ is given, and let $\e > 0$ be given as well.
Let $(\lambda',H') \in \overline{\QQ}$ be another element.
By definition there are sequences
$(\phi_{H_i},H_i)$ and $(\phi_{H'_i}, H'_i)$ converging to
$(\lambda,H)$ and $(\lambda',H')$ respectively.
We have
$$
\align \|\text{Tan}(\phi_{H_i}) - \text{Tan}(\phi_{H'_i}) \|
& = \| H_i \circ \phi_{H_i} - H'_i \circ \phi_{H'_i} \| \\
& \leq \| H_i \circ \phi_{H_i} - H_i \circ \phi_{H'_i} \|
+ \| H_i \circ \phi_{H'_i} - H'_i \circ \phi_{H'_i} \| \\
& = \| H_i \circ \phi_{H_i} - H_i \circ \phi_{H'_i} \| + \| H_i - H'_i \| .
\endalign
$$
By Lemma 3.21, we can find $0 < \delta < \frac{\e}{2}$
and $i_0$ only depending on
the sequence $H_i$ such that: if $\| H_i - H'_i \| < \delta$ and
$d_{C^0} (\phi_{H_i}, \phi_{H'_i})< \delta$ for sufficiently
large $i$, say $i \geq N$, then
$$
\|\text{Tan}(\phi_{H_i}) - \text{Tan}(\phi_{H'_i}) \| \le
\| H_i \circ \phi_{H_i} - H_i \circ \phi_{H'_i} \|  + \| H_i - H'_i \| < \frac{\e}{2} + \frac{\e}{2} = \e
$$
for all $i \geq \max\{i_0,N\}$.
By taking the limit as $i \to \infty$,
this implies
$$\| \overline{\text{Tan}^Q} (\lambda,H) - \overline{\text{Tan}^Q} (\lambda',H') \|
< \e \quad \text{when } \overline{d} (\lambda,\mu) + \| H - H' \| < \delta, $$
proving that $\overline{\text{Tan}^Q}$ is continuous at $(\lambda,H)$.
\qed \enddemo

The images of $\overline{\text{\rm Tan}^Q}$ and
$\overline{\text{\rm Dev}^Q}$ contain $C^\infty_m([0,1] \times
M,\R)$. This is because for any given $F\in C^\infty_m([0,1]
\times M,\R)$, we have the formula
$$
F = \text{Dev}(\phi_F)= - \text{ Tan}(\phi_F^{-1}) \tag 3.23
$$
by (3.6). In fact we will see in Theorem 4.1 that $\text{Im} \,
\overline{\text{\rm Dev}^Q}$ and $\text{Im} \, \overline{\text{\rm
Tan}^Q}$ both contain $C^{1,1}([0,1]\times M,\R)$. We do not know
whether the images of the maps
$$
\overline{\text{Tan}^Q}, \, \overline{\text{Dev}^Q}:
\overline{\QQ} \to L_m^{(1,\infty)}([0,1]\times M,\R) $$
contain the whole $C^0_m([0,1] \times M, \R)$.
Some of these questions will be studied in [Oh7].

\medskip

The power of our definition of the Hamiltonian topology using the
sets (3.7) manifests itself in the proof of the following
theorem.

\proclaim{Theorem 3.23} The set $\overline{\QQ}$ forms a topological group.
\endproclaim

\demo{Proof} We first have to show that composition and inverses on $\overline{\QQ}$ are defined.
The other group properties will follow immediately.
We then show that composition and inverse operation are continuous.

Let $(\lambda,H)$ and $(\mu,F) \in \overline{\QQ}$.
By definition there are sequences $(\phi_{H_i},H_i)$ and $(\phi_{F_i},F_i)$ converging to
$(\lambda,H)$ and $(\mu,F)$ respectively in the Hamiltonian topology.
In particular,
\roster
\item both satisfy
$$
\| H - H_i\|, \quad \| F - F_i\| \to 0
\quad \text{as $i \to \infty$}, \tag 3.24
$$
\item and
$$ \overline d (\lambda, \phi_{H_i}) \to 0, \quad \overline
d (\mu, \phi_{F_i}) \to 0
\quad \text{as $i \to \infty$} . \tag 3.25 $$
\endroster
We know by our earlier remark about $\overline{d}$ that
$$ \overline{d} (\lambda \mu, \phi_{H_i} \phi_{F_i}) \to 0
\quad \text{as } i \to \infty . \tag 3.26 $$
Moreover, we recall
$$
H_i \# F_i = H_i + F_i \circ (\phi_{H_i})^{-1},
$$
and this Hamiltonian generates $\phi_{H_i}\phi_{F_i}$. By
assumption, we have $\|H_i - H\| \to 0$. On the other hand, we
derive
$$
\align
\|F_i \circ (\phi_{H_i})^{-1} - F\circ \lambda^{-1}\|
& \leq \|F_i \circ (\phi_{H_i})^{-1} - F_i \circ \lambda^{-1}\|
+\|F_i \circ \lambda^{-1} - F\circ \lambda^{-1}\| \\
& = \|F_i \circ (\phi_{H_i})^{-1} - F_i \circ \lambda^{-1}\| +
\|F_i - F\|.
\endalign
$$
Here the first term converges to zero by Lemma 3.21
and the second does by assumption. These prove
$$
H_i \# F_i \to H + F\circ \lambda^{-1} \tag 3.27
$$
in the $L^{(1,\infty)}$-topology as $i \to \infty$ under the
assumptions (3.24), (3.25).

Therefore if we define the $L^{(1,\infty)}$-function $H \# F$ by
$$
H \# F := H + F \circ \lambda^{-1},
$$
(3.26) and (3.27) imply that the pair $(\lambda \mu, H \# F)$ is
the limit of the sequence
$$
(\phi_{H_i \# F_i}, H_i \# F_i)
$$
and so lies in $\overline Q$ again. And the above proof also shows
that this limit does not depend on the choices of $H_i, \, F_i$
but depends only on $(\lambda,H)$ and $(\mu, F)$.

Now we define the product of $(\lambda,H)$ and $(\mu, F)$ by
$$
(\lambda,H) \circ (\mu,F) := (\lambda \mu, H \# F). \tag 3.28
$$
When restricted to $\QQ$, this obviously agrees with the usual
definition of composition.

For the inverse, let $(\lambda,H)$ as above.
We know that
$$ \overline{d} \left (\lambda^{-1}, (\phi_{H_i})^{-1} \right )
\to 0 \quad \text{as } i \to \infty. \tag 3.29 $$ Moreover, by the
same proof as for the multiplication, we prove
$$
\lim_{i \to \infty} \overline{H_i} = - H \circ \lambda. \tag 3.30
$$
(One can also prove this by recalling $\overline{H_i} =
-\overline{\operatorname{Tan}^Q}(\phi_{H_i})$ and then using the
continuity of $\overline{\operatorname{Tan}^Q}$ from Proposition
3.22.) Then we define
$$
\overline{H} : = - H \circ \lambda \tag 3.31
$$
which also coincides with the limit (3.30) for any sequence $H_i$
satisfying $\|H - H_i\| \to 0$ and $\overline
d(\lambda,\phi_{H_i}) \to 0$. Now we define the inverse
$$
(\lambda,H)^{-1}: = (\lambda^{-1},\overline H). \tag 3.32
$$
 When
restricted to $\QQ$, this again agrees with the usual definition
of the inverse.

This proves that $\overline{\QQ}$ forms a group under $\circ$,
and it is straightforward to check that all group axioms are satisfied.

\medskip

We now have to show that the group operations in $\overline{\QQ}$ are continuous, i.e.,
that the maps
$$ \matrix
\overline{\QQ} \times \overline{\QQ} & \to & \overline{\QQ}, & \quad
((\lambda,H),(\mu,F)) & \mapsto & (\lambda \mu, H \# F), \\
\overline{\QQ} & \to & \overline{\QQ}, & \quad (\lambda,H) & \mapsto &(\lambda^{-1}, \overline{H})
\endmatrix $$
are continuous with respect to the metric $\overline{d} + \| \cdot \|$.

For the composition, suppose we have two sequences $(\lambda_i,H_i')$
and $(\mu_i,F_i') \in \overline{\QQ}$
converging to $(\lambda,H)$ and $(\mu,F)$ in the metric
$\overline{d} + \| \cdot \|$ on $\overline{\QQ}$ respectively.
We have to show that
$$ \align \overline{d} (\lambda \mu, \lambda_i \mu_i) & \to 0 \quad
\text{as } i \to \infty, \quad \text{and} \\
\| H_i' \# F_i' - H \# F \| & \to 0 \quad \text{as } i \to \infty . \endalign $$
The $C^0$-convergence is again immediate.
For the $\| \cdot \|$-convergence, we compute
$$ \align \| H'_i \# F'_i - H \# F \| & =
\| H'_i + F'_i \circ \lambda_i^{-1} - H - F \circ \lambda^{-1} \| \\
& \le \| H'_i - H \| +
\| F'_i \circ \lambda_i^{-1} - F \circ \lambda^{-1} \| \\
& \le \| H'_i - H \| +
\| F'_i \circ \lambda_i^{-1} - F \circ \lambda_i^{-1} \| \\
& \quad \quad + \| F \circ \lambda_i^{-1} - F \circ \lambda^{-1} \| \\
& = \| H'_i - H \| + \| F'_i  - F \| +  \| F \circ \lambda_i^{-1} - F \circ \lambda^{-1} \|.
\endalign
$$
The first two terms converge to zero by assumption. For the third
term, we derive
$$
\align \| F \circ \lambda_i^{-1} - F \circ \lambda^{-1} \| & \le
\| F \circ \lambda_i^{-1} - F_i \circ \lambda_i^{-1} \| +
\| F_i \circ \lambda_i^{-1} - F_i \circ \lambda^{-1} \| \\
& \quad \quad + \| F_i \circ \lambda^{-1} - F \circ \lambda^{-1} \| \\
& = \| F - F_i \| + \| F_i \circ \lambda_i^{-1} - F_i \circ \lambda^{-1} \|
+ \| F_i - F \|,
\endalign
$$
The first and third term converge to zero by assumption,
and the third term by assumption and Lemma 3.21.
That proves continuity of composition.

For the inverse, $\overline{d} (\lambda^{-1}, \lambda_i^{-1}) \to 0$.
Moreover, it is immediate to check that as in the smooth case (3.18) we have
$$ \| \overline{H_i}- \overline{H} \| =
\| \overline{\text{Tan}^Q}(\lambda)-
\overline{\text{Tan}^Q}(\lambda_i) \| \to 0 $$ by continuity of
$\overline{\text{Tan}^Q}$. That completes the proof. \qed
\enddemo

\proclaim{Corollary 3.24} The set $\QQ \subset \overline{\QQ}$ forms a topological subgroup.
\endproclaim

\demo{Proof} $\QQ$ is a topological subspace of $\overline{\QQ}$
by definition of the latter, and the proof of Theorem 3.23 implies
that $\QQ$ is a subgroup. \qed \enddemo

\proclaim{Corollary 3.25} The evaluation map
$$ \overline{ev}_1^Q : \overline{\QQ} \to \HH ameo (M,\omega) $$
is an open map. The set $\HH ameo(M,\omega)$ forms a topological
group under composition. In particular, $Hameo(M,\omega) \subset
Homeo(M)$ forms a subgroup of $Homeo(M)$.
\endproclaim

\demo{Proof} Theorem 3.23 in particular implies that
left multiplication by an element in $\overline{\QQ}$
is a continuous map $\overline{\QQ} \to \overline{\QQ}$.
By definition, the topology on $\HH ameo (M,\omega)$ is the strongest topology on the set
$Hameo (M,\omega)$ such that the above evaluation map $\overline{ev}_1^Q$ is continuous.
The proof of openness of $\overline{ev}_1^Q$ is now the same as the one for $ev_1$ in Corollary 3.12.

The surjective map
$$ \overline{ev}_1^Q : \overline{\QQ} \to \HH ameo (M,\omega) $$
induces a group structure on $\HH ameo (M,\omega)$ in the obvious
way. In fact, composition in this group is just the usual
composition of maps. The map $\overline{ev}_1^Q$ becomes a
homomorphism of (abstract) groups, which is open, continuous, and
surjective. From this it is straightforward to check that $\HH
ameo (M,\omega)$ indeed forms a topological group.

Since as sets $Hameo(M,\omega)$ coincides with $\HH ameo (M,\omega)$,
$Hameo(M,\omega)$ forms a group as well.
It is immediate that $Hameo(M,\omega)$ with this group structure
forms a subgroup of $Homeo(M)$. \qed \enddemo

We now define the notion of topological Hamiltonian fiber bundles.

\definition{Definition 3.26 [Topological Hamiltonian bundle]}
We call a topological fiber bundle $P \to B$ with fiber $(M,\omega)$ a
topological Hamiltonian bundle
if its structure group can be reduced to the group $Hameo(M,\omega)$.
More precisely, $P\to B$ is a topological Hamiltonian bundle if it allows
a trivializing chart $\{(U_\alpha, \Phi_\alpha)\}$ such that
its transition maps are contained in $Hameo(M,\omega)$.
\enddefinition

Recall that in the smooth case, this definition coincides with
that of a symplectic fiber bundle that carries a fiber-compatible
{\it closed} two form (see [GLS]). It seems to be a very
interesting problem to formulate the corresponding $C^0$-analog to the latter.
We hope to study this issue among others elsewhere.

\definition{Remark 3.27 [Weak Hamiltonian Topology]}
We can define the notion of {\it weak Hamiltonian topology}
similarly to the (strong) Hamiltonian topology. In the sets (3.7),
we just replace the $C^0$-distance of the whole paths by the
$C^0$-distance of the time-one maps only. So in the weak
Hamiltonian topology, we do not have any control over the
$C^0$-convergence of the whole paths other than the time-one maps.
Although this seems natural in light of Proposition 3.6, it turns
out that the weak Hamiltonian topology does not behave as nicely
as the strong Hamiltonian topology. For example, it is unlikely
that the map Tan is continuous with respect to the weak
Hamiltonian topology, or that the sets $\overline{\QQ_{w}}$ and
therefore $Hameo^w(M,\omega)$ defined in the same way as in the
strong case form groups. One can easily verify that Remark 3.8,
Proposition 3.10, Proposition 3.11, Corollary 3.12, and Theorem
4.1 still hold respectively in the weak case, while in Theorem
3.13 only path-connectedness, but not local path-connectedness,
still holds. It seems unlikely that the analog to Theorem 4.5
below holds as well. The strong Hamiltonian topology is obviously
stronger than the weak one, but it is an open question whether
they are indeed different in general.

\enddefinition

\head{\bf \S 4. Basic properties of the group of Hamiltonian homeomorphisms}
\endhead

In this section, we extract some basic properties of the group $Hameo(M,\omega)$
that immediately arise from its definition.
We first note that
$$
Ham(M,\omega) \subset Hameo(M,\omega) \subset Sympeo(M,\omega) \tag
4.1
$$
from their definitions.
The following theorem proves that $Hameo(M,\omega)$ contains all
expected $C^k$-Hamiltonian diffeomorphisms with $k \geq 2$.

\proclaim{Theorem 4.1} The group $Hameo(M,\omega)$ contains all
$C^{1,1}$-Hamiltonian diffeomorphisms. More precisely, if $\phi$
is the time-one map of Hamilton's equation $\dot x = X_H(t,x)$ for a
$C^1$-function $H: [0,1] \times M \to \R$ such that \roster
\item $\|H_t\|_{C^{1,1}} \leq C$,
where $C > 0$ is independent of $t \in [0,1]$, and
\smallskip

\item the map $(t,x) \mapsto dH_t(x), \,[0,1] \times M \to T^*M$
is continuous,
\endroster
then $\phi \in Hameo(M,\omega)$.
\endproclaim
\demo{Proof} Note that any such $C^{1,1}$-function can be
approximated by a sequence of smooth functions $H_i: [0,1] \times M
\to \R$ so that
$$
\|H - H_i\| \to 0, \tag 4.2
$$
where $\|\cdot \|$ denotes the $L^{(1,\infty)}$-norm as before.
On the other hand, the vector fields $X_{H_i}(t,x)$ converge to
$X_H (t,x)$ in $C^{0,1}(TM)$ uniformly over $t \in [0,1]$. Therefore the
flow $\phi_{H_i}^t \to \phi_H^t$ and so $\phi_{H_i}^1 \to
\phi_H^1$ in the $C^0$-topology by the standard existence and continuity
theorem of ODE for Lipschitz vector fields. In particular, this
$C^0$-convergence together with (4.2) implies that the sequence
$(\phi_{H_i},H_i)$ is a Cauchy sequence in $\QQ$ with
$$
\lim_{C^0}\phi_{H_i}^1 = \phi_H^1 = \phi.
$$
Therefore $\phi \in Hameo(M,\omega)$. \qed\enddemo

The following provides an example of an area-preserving homeomorphisms on
a surface that is not $C^1$, but still a Hamiltonian homeomorphism.
Therefore we have the following {\it proper} inclusion relation
$$
Ham(M,\omega) \subsetneq Hameo(M,\omega) \subset Sympeo(M, \omega).
$$
\definition{Example 4.2} We will construct an area-preserving
homeomorphism on the unit disc $D^2$ that is the identity near the
boundary $\partial D^2$ and continuous but not differentiable. By
extending the homeomorphism by the identity on $\Sigma = D^2 \cup
\Sigma \setminus D^2$ to the outside of the disc, we can construct
a similar example on a general surface $\Sigma$ (for example by
choosing $D$ inside the domain of a Darboux chart). Similarly one
can construct such an example in higher dimensions. Furthermore a
slight modification of an example like this combined with
Polterovich's theorem on $S^2$ [P2] provides a sequence $\phi_i$
of Hamiltonian diffeomorphisms on $S^2$ such that $\phi_i \to id$
uniformly but $\|\phi_i\| \to \infty$, which demonstrates that the
Hofer norm function $\phi \mapsto \|\phi\|$ is not continuous in
the $C^0$-topology on $Ham(M,\omega)$.

Let $(r,\theta)$ be polar coordinates on $D^2$. Then the standard
area form is given by
$$
\Omega = r\, dr \wedge d\theta.
$$
Consider maps $: D^2 \to D^2$ of the form
$$
\phi_{\rho}: (r,\theta) \mapsto (r, \theta + \rho(r)),
$$
where $\rho: (0,1] \to [0, \infty)$ is a smooth function that satisfies
for some small $\e > 0$
\medskip
\roster
\item $\rho' < 0$ on $(0,1-\e)$, $\rho \equiv 0$ on $[1-\e,1]$, and
\item $\lim_{r \to 0^+} r \rho'(r) = -\infty$.
\endroster
\medskip

It follows that $\phi_\rho$ is smooth except at the origin at
which $\phi_\rho$ is continuous but not differentiable. Obviously
the map $\phi_{-\rho}$ is the inverse of $\phi_\rho$ which shows
that it is a homeomorphism. Furthermore we have
$$
\phi_\rho^*(r\, dr\wedge d\theta) = r\,dr\wedge d\theta \quad \text{on } \,
D^2 \setminus \{0\},
$$
which implies that $\phi_\rho$ is area-preserving.

Now it remains to show that if we choose $\rho$ suitably,
$\phi_\rho$ becomes a Hamiltonian homeomorphism.
We will in fact consider {\it time-independent} Hamiltonians for this purpose.
Consider the isotopy
$$ t\in [0,1] \mapsto \phi_{t\rho} \in Homeo^\Omega(D^2). $$
A straightforward calculation shows that a corresponding (not necessarily normalized)
Hamiltonian is given by the time-independent function
$$ H_\rho(r,\theta) = - \int_1^r s \rho(s)\, ds. $$
The $L^{(1,\infty)}$-norm of $H_\rho$ becomes
$$ \int_0^1 s \rho(s)\, ds . $$
Choose any $\rho$ so that the integral becomes finite, e.g.
$\rho(r) = \frac{1}{\sqrt{r}}$ near $r = 0$. Now we choose any
smoothing sequence $\rho_n$ of $\rho$ by regularizing $\rho$ at
$0$, and consider the corresponding Hamiltonians $H_{\rho_n}$ and
their time one-maps $\phi_{\rho_n}$. Then it follows that
$(\phi_{H_{\rho_n}},H_{\rho_n})$ is a Cauchy sequence in the
Hamiltonian topology and $\phi_{\rho_n} \to \phi_{\rho}$ in the
$C^0$-topology. So $\phi_\rho$ is a Hamiltonian
homeomorphism that is neither differentiable nor Lipschitz at $0$.
\enddefinition

The following question seems to be one of fundamental importance
(See Conjectures 5.3 and 5.4 later).

\definition{Question 4.3} In Example 4.2, consider $\rho$ such that
$$
\int_{0^+}^1 s \rho(s)\, ds = +\infty.
$$
Is the homeomorphism $\phi_\rho$ still contained in $Hameo(M,\omega)$?
\enddefinition

The following theorem is the $C^0$-version of the
well-known fact that $Ham(M,\omega)$ is a normal subgroup of $Symp_0(M,\omega)$.

\proclaim{Theorem 4.4} $Hameo(M,\omega)$ is a normal subgroup of $Sympeo(M,\omega)$.
\endproclaim
\demo{Proof} We have to show
$$
\psi h \psi^{-1} \in Hameo(M,\omega)
$$
for any $h \in Hameo(M,\omega)$ and $\psi \in Sympeo(M,\omega)$.
By definition, there are sequences $(\phi_{H_i},H_i) \in \QQ$ and
$\psi_i \in Symp(M,\omega)$ such that
$$
h = \lim_{C^0}\phi_{H_i}^1 \quad \text{and} \quad \lim_{C^0}\psi_i = \psi .
$$
Let $\phi_i = \phi_{H_i}^1$.
Recall from (3.3) that $\psi_i^{-1} \phi_i \psi_i$ is generated by $H_i \circ \psi_i$ for all $i$.
It therefore suffices to prove that $(\psi_i^{-1} \phi_i \psi_i, H_i \circ \psi_i)$
is a Cauchy sequence in $\QQ$ and
$\lim_{C^0} \psi_i^{-1} \phi_i \psi_i = \psi^{-1} h \psi$.
The $C^0$-convergence of the paths and time-one maps is obvious.
Hence it remains to prove that
$H_i \circ \psi_i$ is a Cauchy sequence in the $L^{(1,\infty)}$-topology,
$$
\|H_i \circ \psi_i - H_j \circ \psi_j\| \to 0 \quad \text{as }\, i, \, j \to \infty. \tag 4.4
$$
But
$$
\| H_i \circ \psi_i - H_j \circ \psi_j \|
\leq \|H_i\circ \psi_i - H_j\circ \psi_i \| + \| H_j\circ \psi_i - H_j \circ \psi_j \| \to 0.
$$
Here the first term goes to zero as $\|H_i\circ \psi_i - H_j\circ \psi_i \|
=\|H_i - H_j\| \to 0$ by assumption,
and the second does by assumption and by Lemma 3.21(2)
(by viewing the $\psi_i$ as constant paths).
That finishes the proof. \qed\enddemo

The following is an important property of $Hameo(M,\omega)$,
which demonstrates that it is the `correct' $C^0$-counterpart of $Ham(M,\omega)$.

\proclaim{Theorem 4.5} $\HH ameo(M,\omega)$ is path-connected and locally path-connected.
Consequently, $Hameo(M,\omega)$ is path-connected and we have
$$
Hameo(M,\omega) \subset Sympeo_0(M,\omega) \subset
Sympeo(M,\omega)\cap Homeo_0^\Omega(M).
$$
\endproclaim
\demo{Proof}
Let $h \in \HH ameo(M,\omega)$.
For the path-connectedness of $\HH ameo(M,\omega)$,
it suffices to prove that $h$ can be connected to the identity by
a Hamiltonian continuous path
$\ell : [0,1] \to \HH ameo(M,\omega)$ such that
$\ell(0) = id$ and $\ell(1) = h$.

By definition, there exists a sequence $(\phi_{H_i},H_i) \in \QQ$
converging to an element $(\lambda,H) \in \overline{\QQ}$, and $h
= \overline{ev}_1^Q (\lambda,H) = \lambda(1) = \lim_{C^0}
\phi_{H_i}^1$. As in Theorem 3.13 consider the  $H_i^s$ generating
the Hamiltonian paths $t \mapsto \phi_{H_i^s}^t = \phi_{H_i}^{st}$
for all $s \in [0,1]$ and all $i$. By the same arguments as in
Theorem 3.13 we have
$$ \align
\overline d(\phi_{H_i^s},\phi_{H_{i'}^s}) \leq \overline d(\phi_{H_i},\phi_{H_{i'}})
& \to 0 \quad \text{as } i, \, i' \to \infty, \quad \text{and} \\
\| H_i^s - H_{i'}^s \| \leq \| H_i - H_{i'} \| & \to 0 \quad
\text{as } i, \, i' \to \infty .
\endalign $$
So $(\phi_{H_i^s}, H_i^s)$ is a Cauchy sequence in the Hamiltonian
topology. Denote by $(\lambda^s,H^s) \in \overline{\QQ}$ its
limit, and note that $\lambda^s$ is nothing but the path $t
\mapsto \lambda (st)$. By the above, $\ell(s) = \overline{ev}_1^Q
(\lambda^s, H^s) = \lambda(s) \in \HH ameo(M,\omega)$ for all $s
\in [0,1]$, and $\ell(0) = id$, $\ell(1) = h$. It remains to show
that $\ell$ is continuous with respect to the Hamiltonian topology
on $\HH ameo(M,\omega)$.

Now $\ell$ factors through
$$ [0,1] \to \overline{\QQ} \to \HH ameo(M,\omega),
\quad s \mapsto (\lambda^s,H^s) \mapsto \overline{ev}_1^Q (\lambda^s, H^s) = \ell (s). $$
By definition of the topology on $\HH ameo(M,\omega)$
it suffices to show that the first map is continuous, that is,
that $s \mapsto (\lambda^s,H^s)$ is continuous with respect to the standard metric on $[0,1]$
and the product metric $\overline{d} + \| \cdot \|$ on $\overline{\QQ}$.
But
$$
\align \overline{d} \Big( (\lambda^s,H^s), (\lambda^{s'},H^{s'}) \Big) & =
\| H^s - H^{s'} \| + \overline{d} ( \lambda^s, \lambda^{s'} ) \\
& = \lim_{i \to \infty}  \| H_i^s - H_i^{s'} \|  + \max_{t \in
[0,1]} \overline d (\lambda (st), \lambda (s't)) .
\endalign $$
Let $\e > 0$. Note that if we consider the functions
$\zeta_1(t) = ts$ and $\zeta_2(t)= ts'$, we see that
$$
\|\zeta_1 - \zeta_2\|_{ham} = 2 |s-s'|.
$$
Therefore it follows from Lemma 3.21 that we can find $\delta > 0$ and $i_0$ sufficiently large such that
$$\| H_i^s - H_i^{s'} \| < \frac{\e}{2} $$
when $| s - s' | < \delta$ and $i \geq i_0$,
and therefore
$$ \lim_{i \to \infty} \| H_i^s - H_i^{s'} \|  < \frac{\e}{2} $$
when $| s - s' | < \delta$. For the second term, use continuity of
$\lambda$ and $\lambda^{-1}$ to see that by making $\delta$
smaller if necessary,
$$ \overline d ( \lambda (st), \lambda (s't) ) < \frac{\e}{2} $$
when $| st - s't | \leq | s - s' | < \delta$.
That proves continuity of $\ell$,
and hence completes the proof of path-connectedness of $\HH ameo(M,\omega)$.

For the proof of local path-connectedness, we can, using Corollary 3.25,
combine the above proof with the ideas in the proof of Theorem 3.13.
Since the proof is essentially the same, we leave the details to the reader.

Now as sets, $Hameo(M,\omega)$ coincides with $\HH ameo(M,\omega)$.
Note that the path $\ell$ constructed above is a topological Hamiltonian path.
Since a topological Hamiltonian path is in particular a continuous path with respect to the $C^0$-topology,
this implies path-connectedness of $Hameo(M,\omega)$.
The other statements about $Hameo(M,\omega)$ follow from this immediately.
That completes the proof.
\qed\enddemo

It follows immediately from the $L^{(1,\infty)}$-Approximation Lemma (Appendix 2)
that given any Cauchy sequence in $\QQ$, we may assume that each path in the sequence is boundary flat.
This implies that the concatenation of two topological Hamiltonian path is again a topological Hamiltonian path.
So in fact we have proved that $Hameo (M,\omega)$ is path-connected by {\it topological Hamiltonian path}.
The proof involves the boundary flattening procedure,
and therefore only works in the $L^{(1,\infty)}$-topology and not in the $L^\infty$-topology.
As remarked above, this is one indication that the $L^{(1,\infty)}$-topology, not the $L^\infty$-topology,
is the correct topology for the study of $C^0$-Hamiltonian geometry.

\definition{Question 4.6}
Is $Hameo(M,\omega)$ locally path-connected?
\enddefinition

\medskip

Recall that by (4.1) we have $ Hameo(M,\omega) \subset Sympeo(M,\omega)$.
But note that a priori it is not obvious whether $Hameo(M,\omega)$
is different from $Sympeo(M,\omega)$.
In fact, if one naively takes just the $C^0$-closure of $Ham(M,\omega)$,
then it can end up becoming the whole $Sympeo(M,\omega)$.
We refer to [Bt] for a nice observation that this is really the case for $Ham^c(\R^{2n})$.
We refer to section 6 for further discussion on this phenomenon.

In the next section, we will study the case $\dim M = 2$.
Here we want to state the following theorem which is an immediate
application of Arnold's conjecture.

\proclaim{Theorem 4.7} Let $(M,\omega)$ be a closed symplectic
manifold. Then any $C^0$-limit of Hamiltonian diffeomorphism  has
a fixed point. In particular, any Hamiltonian homeomorphism has a
fixed point.
\endproclaim
\demo{Proof} Let $h = \lim_{C^0} \phi_i$ for a sequence $\phi_i
\in Ham(M,\omega)$. We prove the theorem by contradiction.
Suppose $h$ has no fixed point. Denote
$$
d^h_{min}: = \inf_{x \in M} d(x,h(x)).
$$
By compactness of $M$ and since $h$ has no fixed point, $d^h_{min}
> 0$. But each $\phi_i$ must have a fixed point $x_i$ by the
Arnold Conjecture, which was proven in [FOn], [LT] or [Ru]. Hence
$$
\overline d (h, \phi_i) \geq d (h(x_i), \phi_i(x_i)) = d
(h(x_i), x_i) \geq d^h_{min}>0
$$
for all $i$. On the other hand, we have
$$
\lim_{i \to \infty} \overline d (h, \phi_i) = 0 ,
$$
which gives rise to a contradiction. \qed
\enddemo

\proclaim{Corollary 4.8} Suppose that $(M,\omega)$ carries a
symplectic diffeomorphism $\psi \in Symp_0(M,\omega)$ (or
equivalently, $\psi \in Sympeo_0(M,\omega)$) that has no fixed
point. Then $\psi \not \in Hameo(M,\omega)$, and in particular we
have
$$
Hameo(M,\omega) \subsetneq Sympeo_0(M,\omega).
$$
\endproclaim

An example of a symplectic manifold $(M,\omega)$ satisfying the hypothesis of Corollary 4.8
is the torus $T^{2n}$ with the standard symplectic form $\omega_0$.
Recall that by identifying $\alpha \in T^{2n}$ with the rotation $x \mapsto x + \alpha$,
we can identify $T^{2n}$ with a subgroup of $Symp_0(T^{2n},\omega_0)$,
$$ T^{2n} \hookrightarrow Symp_0(T^{2n},\omega_0) . $$
By Theorem 4.7, we have
$$ T^{2n} \cap Hameo(T^{2n},\omega_0) = \{ id \} . $$
It follows that $Hameo(T^{2n},\omega_0) \subsetneq Sympeo_0(T^{2n},\omega_0)$.

\head{\bf \S 5. The two dimensional case}
\endhead

In this section, we will mainly study the case $\dim M = 2$.
The first question would be what the relation between the group
$Homeo^\Omega(M)$ ($Homeo^\Omega_0(M)$) and its subgroup
$Sympeo(M,\omega)$ ($Sympeo_0(M,\omega)$) is.
By definition of $Sympeo(M,\omega)$, in two dimensions this question
boils down to the approximability of area-preserving homeomorphisms
by area-preserving diffeomorphisms. We refer readers to [Oh6] for the
precise statements and proofs but state their consequence here
because our discussion in this section will be based on this theorem.

\proclaim{Theorem 5.1 [Oh6]} Let $M$ be a compact orientable surface
without boundary and $\omega = \Omega$ be an area form on it. Then we have
$$
Sympeo(M,\omega) = Homeo^\Omega(M), \quad Sympeo_0(M,\omega)
= Homeo^\Omega_0(M).
$$
\endproclaim

Next we study the relation between $Hameo(M,\omega)$ and
$Sympeo_0(M,\omega) = Homeo_0^\Omega(M)$. We will prove
that the subgroup $Hameo(M,\omega) \subset Sympeo_0(M,\omega)$ is
indeed a proper subgroup if $M \neq S^2$.
The proof will use the  mass flow homomorphism
for area-preserving homeomorphisms on a surface, which we
recalled in section 2 in the general context of measure-preserving
homeomorphisms. The mass flow homomorphisms can be defined for any
isotopy of measure-preserving homeomorphisms preserving a {\it good}
measure, e.g., the Liouville measure on a symplectic manifold $(M,\omega)$.
The mass flow homomorphism reduces to the dual
version of the flux homomorphism for volume-preserving {\it
diffeomorphisms} on a {\it smooth} manifold [T]. Of course in two
dimensions, the flux homomorphism coincides with the symplectic
flux homomorphism, and so we can compare the mass flow homomorphism
and the symplectic flux. One crucial point of considering the mass
flow homomorphism instead of the flux homomorphism
is that it is defined for an isotopy of area-preserving
{\it homeomorphisms}, not just for diffeomorphisms.

We first recall the definition of the symplectic flux homomorphism.
Denote by
$$
\PP(Symp_0(M,\omega),id)
$$
the space of {\it smooth} paths $c: [0,1] \to Symp_0(M,\omega)$
with $c(0) = id$. This naturally forms a group. For each given
$c\in \PP(Symp_0(M,\omega),id)$, the Flux of $c$ is defined by
$$
\PP(Symp_0(M,\omega),id) \to H^1(M,\R), \quad
Flux(c)= \int_0^1\dot c\, \rfloor \omega \, dt . \tag 5.1
$$
This depends only on the homotopy class, relative to the end
points, of the path $c$ and therefore projects down to the
universal covering space
$$
\pi_\omega: \widetilde{Symp_0}(M,\omega) \to Symp_0(M,\omega), \quad [c]
\mapsto c(1), \tag 5.2
$$
where
$$
\widetilde{Symp_0}(M,\omega): =\{\, [c] \, \mid \, c
\in\PP(Symp_0(M,\omega),id)\}.
$$
Here $[c]$ is the homotopy class of $c$ relative to fixed end points.
We recall that $Symp_0(M,\omega)$ is locally contractible [W]
and so $\widetilde{Symp}_0(M,\omega)$ is indeed the universal covering
space of $Symp_0(M,\omega)$.
If we put
$$ \Gamma_\omega = \text{Flux} \, \left( \ker \left(
\pi_\omega \colon \widetilde{Symp_0}(M,\omega) \to Symp_0(M,\omega)\right) \right) , $$
we obtain by passing to the quotient the group homomorphism
$$ \text{flux} \, \colon Symp_0(M,\omega) \to H^1(M,\R) / \Gamma_\omega . \tag 5.3 $$
The maps (5.1) and (5.3) are also known to be surjective [Ba].

It is also shown in [Fa, Appendix A.5] that $\text{Flux} (c) \in H^1(M,\R)$ is the Poincar\'e dual to
the mass flow homomorphism $\widetilde\theta (c) \in H_1(M,\R)$ recalled in section 2
(after normalizing $\omega$ so that $\int_M \omega = 1$).
Since it is also well-known [Ba] that
$$
\align \widetilde{Ham}(M,\omega) & = \ker \text{Flux}, \\
Ham(M,\omega) & = \ker \text{flux},
\endalign
$$
we derive
$$
Ham(M,\omega) \subset \ker \theta \cap Symp_0(M,\omega). \tag 5.4
$$

\proclaim{Theorem 5.2} Let $(M,\omega)$ be a closed orientable surface,
where $\omega = \Omega$ is a symplectic (or area) form on $M$.
Then we have
$$
Hameo(M,\omega) \subset \ker \theta \cap Sympeo_0(M,\omega). \tag 5.5
$$
In particular, if $M \neq S^2$, we have
$$
Hameo(M,\omega) \subsetneq Sympeo_0(M,\omega) = Homeo^\Omega_0(M). \tag 5.6
$$
\endproclaim
\demo{Proof} Recall (4.1) that $Hameo(M,\omega) \subset
Sympeo_0(M)$. On the other hand, (5.4) implies
$\theta|_{Ham(M,\omega)}\equiv 0$. From continuity of $\theta$
(Theorem 2.2) and the definition of $Hameo(M,\omega)$ we derive
$\theta|_{Hameo(M,\omega)} \equiv 0$. That proves (5.5).

By the surjectivity of the Flux, the map
$\theta|_{Sympeo_0(M,\omega)} : Sympeo_0(M,\omega) \to H_1(M,\R) / \Gamma$ is surjective.
So $\ker \theta|_{Sympeo_0(M,\omega)} \subsetneq Sympeo_0(M,\omega)$ when
$H_1(M,\R) \neq 0$ (and therefore $H_1(M,\R) / \Gamma \not= 0$ since $\Gamma$ is discrete)
which is the case for $M \neq S^2$.
That proves the last statement.
\qed\enddemo

This theorem verifies that $Hameo(M,\omega)$ is a {\it proper}
normal subgroup of $Sympeo(M,\omega)$, at least in two dimensions
if $M\neq S^2$.

We now propose the following conjecture

\proclaim{Conjecture 5.3} $Hameo(M,\omega)$ is a proper subgroup of $\ker \theta$ in general.
In particular for $M = S^2$ with $\Omega = \omega$, $Hameo(S^2,\omega)$ is a
proper normal subgroup of $Sympeo_0(S^2,\omega) = Homeo^\Omega_0(S^2)$.
\endproclaim

The affirmative answer to this conjecture will answer to Question 2.3 negatively
and settle the simpleness question of $Homeo^\Omega_0(S^2)$, which has been open
since Fathi's paper [Fa] appeared. In fact, this conjecture is
an immediate corollary of the following more concrete conjecture

\proclaim{Conjecture 5.4} The answer to Question 4.3 on $S^2$ is negative,
at least for a suitable choice of $\rho$.
\endproclaim

The results of this section can be generalized to higher dimensions in many cases.
We first recall the flux homomorphism for volume-preserving
{\it diffeomorphisms} on a {\it smooth} manifold [T].
Let $\Omega$ be a volume form on $M$ and denote by
$$
\PP(Diff_0^\Omega(M),id)
$$
the space of smooth paths $c: [0,1] \to Diff_0^\Omega(M)$, the group of diffeomorphisms
preserving the volume form $\Omega$, with $c(0) = id$.
This also naturally forms a group. For each given $c\in
\PP(Diff_0^\Omega(M),id)$, the Volume Flux of $c$ is defined by
$$
\PP(Diff_0^\Omega(M),id) \to H^{2n-1}(M,\R), \quad
\tilde{V} (c) = \int_0^1\dot c\, \rfloor \Omega \, dt .
$$
This depends only on the homotopy class relative to the end points
of the path $c$ and therefore projects down to the universal
covering space
$$
\pi_\Omega: \widetilde{Diff_0^\Omega}(M) \to Diff_0^\Omega(M), \quad [c] \mapsto c(1) ,
$$
where
$$
\widetilde{Diff_0^\Omega}(M): =\{\, [c] \, \mid \, c \in\PP(Diff_0^\Omega(M,\omega),id)\}.
$$
Here $[c]$ again denotes the homotopy class of $c$ relative to fixed end points.
It is well-known that $Diff_0^\Omega(M)$ is locally contractible and
so $\widetilde{Diff_0^\Omega}(M)$ is indeed the universal covering
space of $Diff_0^\Omega(M)$.
If we put
$$ \Gamma_\Omega = \widetilde{V} \, \left( \ker \left(
\pi_\Omega \colon \widetilde{Diff_0^\Omega}(M) \to Diff_0^\Omega(M)\right) \right) , $$
we obtain by passing to the quotient the group homomorphism
$$ V \, \colon Diff_0^\Omega(M) \to H^{2n-1}(M,\R) / \Gamma_\Omega , $$
to which we also refer to as the {\it (volume) flux homomorphism}.

In fact [Fa], $\tilde{V} (c) \in H^{2n-1}(M,\R)$ is the Poincar\'e dual to
the mass flow homomorphism $\widetilde\theta (c) \in H_1(M,\R)$
(after normalizing $\Omega$ so that $\int_M \Omega = 1$).

Now let $\Omega = \frac{1}{n!} \omega^n$ be the Liouville volume
form. An easy calculation [Ba] shows that
$$
\widetilde{V} (c) = \frac{1}{(n-1)!} \Big( \text{Flux} (c) \Big) \wedge
\omega^{n-1} . \tag 5.7
$$
So (5.4) holds in any dimension,
$$
Ham(M,\omega) \subset \ker \theta \cap Symp_0(M,\omega).
$$
By reexamining the proof of Theorem 5.2, we see that (5.5) holds as well, i.e.,
$$
Hameo(M,\omega) \subset \ker \theta \cap Sympeo_0(M,\omega)
$$
for any closed symplectic manifold $(M,\omega)$. We also see that
$$
Hameo(M,\omega) \subsetneq Sympeo_0(M,\omega)
$$
if $\theta|_{Sympeo_0(M,\omega)} \colon Sympeo_0(M,\omega) \to H_1(M,\omega)/\Gamma$ is nontrivial.
By (5.7) and surjectivity of the Flux, we see that this condition is satisfied if
$$
\wedge \omega^{n-1} \colon H^1(M,\R) \to H^{2n-1}(M,\R) \tag 5.8
$$
is nontrivial.
Since the map (5.8) is easily seen to be surjective, the latter condition is satisfied whenever
$H^{2n-1} (M,\R) \cong H_1 (M,R)$ (by Poincar\'e duality) is nontrivial.
This holds for example for the torus $T^{2n}$ and therefore gives another proof of
$Hameo(T^{2n},\omega_0) \subsetneq Sympeo_0(T^{2n},\omega_0)$,
which was also a consequence of Corollary 4.8.
We summarize these results in the following theorem.

\proclaim{Theorem 5.5} Let $(M,\omega)$ be a closed symplectic manifold.
Then we have
$$
Ham(M,\omega) \subset \ker \theta \cap Symp_0(M,\omega),
$$
and
$$
Hameo(M,\omega) \subset \ker \theta \cap Sympeo_0(M,\omega). \tag 5.9
$$
If in addition
$$
H_1(M,\R) \cong H^{2n-1}(M,\R)
$$
is nontrivial, then
$$ Hameo(M,\omega) \subsetneq Sympeo_0(M,\omega) \subset Homeo^\Omega_0(M). $$
\endproclaim

\medskip
\head \bf \S6. The non-compact case and open problems
\endhead

So far we have assumed that $M$ is closed. In this section, we
will indicate the necessary changes to be made for the open case
where $M$ is either noncompact or with boundary or both.

There are two possible definitions of compactly supported
Hamiltonian diffeomorphisms in the literature. In this paper, we
will treat the more standard version, which we call {\it compactly
supported Hamiltonian diffeomorphisms}.

Here is the definition of compactly supported Hamiltonian
diffeomorphisms which is mostly used in the literature so far. We
denote $Symp^c(M,\omega) \subset Diff^c(M,\omega)$ the set of
compactly supported symplectic diffeomorphisms.

\definition{Definition 6.1}
We say that a smooth path $\lambda:[0,1] \to Symp^c(M,\omega)$ is
a compactly supported Hamiltonian path if $\lambda = \phi_H$ for a
Hamiltonian function $H: [0,1] \times M \to \R$ such that $H$ is
compactly supported in $\operatorname{Int}(M)$ and $\phi =
\phi_H^1$, where $\text{supp}(H)$ is defined by
$$
\text{supp}(H) = \overline{\bigcup_{t \in [0,1]}\text{supp}(H_t)}.
$$
We define
$$
\PP^{ham}(Symp^c(M,\omega), id)
$$
to be the set of such $\lambda$'s. A compactly supported
symplectic diffeomorphism $\phi$ is a {\it compactly supported
Hamiltonian diffeomorphism} if $\phi = ev_1(\lambda)$ for a
$\lambda \in \PP^{ham}(Symp^c(M,\omega), id)$. We denote
$$
Ham^c(M,\omega)= ev_1(\PP^{ham}(Symp^c(M,\omega), id)).
$$
\enddefinition

We now give a description of the Hamiltonian topology on
$\PP^{ham}(Symp^c(M,\omega), id)$ and $Ham^c(M,\omega)$.

Let $K \subset \operatorname{Int}(M)$ be a compact subset. We
denote by $Symp_K(M,\omega)$ to be the subset of
$Symp^c(M,\omega)$ and then by definition
$$
Symp^c(M,\omega) = \bigcup_{K \subset \operatorname{Int}M;compact}
Symp_K(M,\omega).
$$
We denote by
$$
\PP^{ham}(Symp_K(M,\omega),id)
$$
the set of $\lambda \in \PP^{ham}(Symp^c(M,\omega),id)$ with
$$
\operatorname{supp}(\lambda(t)) \subset K \quad \text{for all $t
\in [0,1]$}.
$$
The Hamiltonian topology on $\PP^{ham}(Symp_K(M;\omega),id)$ is
equivalent to the metric topology thereon induced by the metric
$$
d_{ham,K}(\lambda_0,\lambda_1) = \overline d(\lambda_0,\lambda_1)
+ \operatorname{leng}(\lambda_0^{-1}\lambda_1)
$$
(Proposition 3.10), where $\overline d$ is the $C^0$-metric on
$\PP(Homeo^c(M),id)$. By definition,
$$
\PP^{ham}(Symp^c(M,\omega),id) = \bigcup_{K\subset
\operatorname{Int}M;compact} \PP^{ham}(Symp_K(M,\omega),id).
$$

We then define $Ham_K(M,\omega)$ to be the image
$$
Ham_K(M,\omega) = ev_1(\PP^{ham}(Symp_K(M,\omega),id)).
$$
\definition{Definition 6.2} Suppose $M$ is either noncompact or
with boundary $\part M \neq \emptyset$. Then \roster \item the
strong Hamiltonian topology $\PP^{ham}(Symp^c(M,\omega),id)$ is
the direct limit topology of the directed system
$$
\{\PP^{ham}(Symp^c(M,\omega),id) \mid K\subset
\operatorname{Int}M, \, \text{compact}\}.
$$
\smallskip

\item We define the {\it Hamiltonian topology} of
$Ham^c(M,\omega)$ to be the strongest topology thereon such that
the evaluation map
$$
ev_1 \colon \PP^{ham}(Symp^c(M,\omega),id) \to Symp^c(M,\omega).
$$
is continuous. We denote the resulting topological space by $\HH
am^c(M,\omega)$.
\endroster
\enddefinition

Note that by definition we have
$$
Ham^c(M,\omega) = \bigcup_{K \subset \operatorname{Int}M;compact}
Ham_K(M,\omega).
$$
An easy exercise, using the commutative diagram
$$
\matrix ev_1 & \colon \PP^{ham}(Symp_K(M,\omega),id)
&\longrightarrow &Symp_K(M,\omega) \\
& \downarrow &\quad & \downarrow \\
ev_1 & \colon  \PP^{ham}(Symp^c(M,\omega),id) & \longrightarrow &
Symp^c(M,\omega),
\endmatrix
$$
shows that the Hamiltonian topology of $Ham^c(M,\omega)$ is
equivalent to the direct limit topology of $Ham_K(M,\omega)$ over
$K$.

Now the developing map $\operatorname{Dev}$ has the form
$$
\operatorname{Dev} : \PP^{ham}(Symp^c(M,\omega),id) \to
C_c^\infty([0,1] \times M,\R).
$$
Here $C_c^\infty([0,1] \times M,\R)$
is the set of smooth functions such that
$$
\overline{\cup_{t \in [0,1]}\operatorname{supp}(H_t)} \subset
\operatorname{Int}(M)
$$
is compact.

We also consider the inclusion map
$$
\align \iota_{ham}: \PP^{ham}(Symp^c(M,\omega),id) & \to
\PP(Symp^c(M,\omega),id) \\
& \to   \PP(Homeo^c(M),id).
\endalign
$$
The unfolding map $(\iota_{ham},Dev)$ has the image
$$
\QQ := \operatorname{Image} (\iota_{ham},Dev) \subset
\PP(Homeo^c(M),id) \times L^{(1,\infty)}_c([0,1] \times M,\R),
$$
Similarly we define
$$
\QQ_K := \operatorname{Image} (\iota_{ham,K},Dev_K) \subset
\PP(Homeo_K(M),id) \times L^{(1,\infty)}_K([0,1] \times M,\R)
$$
which has the unique topology induced by the metric topology on
$\QQ_K$. Now we equip $\QQ$ the direct limit topology of $\QQ_K$.
Then it follows that the unfolding map canonically extends to the
union
$$
\overline{\QQ}:= \bigcup_{K\subset \operatorname{Int}M;compact}
\overline\QQ_K
$$
in that we have the following continuous projections
$$
\align \overline{\iota}_{ham}^Q & : \overline{\QQ} \to
\PP(Homeo^c(M),id) \tag 6.2\\
\overline{\operatorname{Dev}^Q} & : \overline{\QQ} \to
L^{(1,\infty)}_c([0,1]\times M,\R) \tag 6.3
\endalign
$$
with respect to the direct limit topology of $\overline\QQ$ and
the similar topology on the targets. We would like to remark that
$\overline \QQ$ is {\it not} the closure of the metric topology on
$\PP(Homeo^c(M),id) \times L^{(1,\infty)}_c([0,1] \times M,\R)$ :
the latter product space is not a complete metric space.

By definition we have the extension of the evaluation map
$$
ev_1:\PP^{ham}(Symp^c(M,\omega),id) \to Symp^c(M,\omega) \to
Homeo^c(M)
$$
to
$$
\overline{ev}_1^Q: \overline\QQ \to Homeo^c(M). \tag 6.4
$$

\definition{Definition 6.3}
We define the set
$$
\align \PP^{ham}(Sympeo_K(M,\omega),id) & : =
\overline{\iota}_{ham}^Q(\overline \QQ_K) \subset
\PP(Homeo_K(M),id) \\
\PP^{ham}(Sympeo^c(M,\omega),id) & : = \overline{\iota}_{ham}^Q
(\overline \QQ) \subset \PP(Homeo^c(M),id)
\endalign
$$
and call any element of $\PP^{ham}(Sympeo^c(M,\omega),id)$ a
compactly supported {\it topological Hamiltonian path}. Again we
equip the latter with the direct limit topology of the metric
topologies on $\PP^{ham}(Sympeo_K(M,\omega),id)$. We call this the
Hamiltonian topology on $\PP^{ham}(Sympeo^c(M,\omega),id)$.
\enddefinition

Then the set of compactly supported {\it Hamiltonian
homeomorphisms} is defined by
$$
\aligned Hameo^c(M,\omega) & = \{h \in Homeo(M) \mid h =
\overline{ev}_1(\lambda),
\\
&{ } \qquad \lambda\in \PP^{ham}(Sympeo^c(M,\omega),id)\}
\endaligned
\tag 6.5
$$

\definition{Definition 6.4} We define
$$
Hameo_K(M,\omega) = \overline{ev}_1^Q(\overline\QQ), \quad
(\lambda,H) \to \lambda(1)
$$
and then
$$
Hameo^c(M,\omega) = \bigcup_{K\subset \operatorname{Int}M
;compact} Hameo_K(M,\omega).
$$
We call the {\it Hamiltonian topology} on $Hameo^c(M,\omega)$ the
direct limit topology the metric topologies on
$Hameo_K(M,\omega)$.
\enddefinition
With these definitions, the analogs to all the results stated in
section 2-5 still hold. For example, the following can be proved
in the same way as Theorem 4.4 and Theorem 4.5.

\proclaim{Theorem 6.5} The group $Hameo^c(M,\omega)$ is a
path-connected normal subgroup of $Sympeo_0^c(M,\omega)$.
\endproclaim

We would like to point out that this theorem is a sharp contrast
to the following interesting observation by S. Bates [Bt]: {\it if
one takes just the $C^0$-closure instead, not with respect to the
Hamiltonian topology, of $Ham^c(\R^{2n},\omega_0)$,
$Hameo^c(\R^{2n},\omega)$ is the whole
$Sympeo^c(\R^{2n},\omega_0)$ even if $Symp(\R^{2n},\omega_0)$ has
many connected components}. This is another evidence the
Hamiltonian topology is the right topology to take for the study
of topological Hamiltonian geometry.

In relation to this definition, we would just like to mention one
result by Hofer [H2] on $\R^{2n}$ :
$$
\|\phi^{-1}\psi\| \leq C\, \text{diam}(\text{supp}(\phi^{-1}\psi))
\|\phi^{-1}\psi\|_{C^0}, \tag 6.6
$$
where $C$ is a constant with the bound $C \leq 128$. This in
particular implies that the $C^0$-topology is stronger than the
Hofer topology on $Ham^c(\R^{2n},\omega_0)$ {\it if
$\text{supp}(\phi^{-1}\psi)$ is controlled.}

\medskip

Finally we list the problems which arise immediately from the
various definitions introduced in this paper, and seem to be
interesting to investigate. These will be subjects of future
study.

\definition{Problems}
\roster \item Describe the closed set of length minimizing paths
in terms of the geometry and dynamics of the Hamiltonian flows.
\item Describe the images of $\overline{\text{Tan}^Q},
\overline{\text{Dev}^Q}$ of $\overline{\QQ}$ in
$L^{(1,\infty)}_m([0,1]\times M,\R)$. \item Study the structure of
the flow of Hamiltonian homeomorphisms in terms of the
$C^0$-Hamiltonian dynamical system or as the high dimensional
generalization of area-preserving homeomorphisms with vanishing
mass flow or zero mean rotation vector. \item Does the identity
$[Sympeo_0, Sympeo_0] = Hameo$ hold? Is $Hameo$ simple? \item
Further investigate the above Hofer's inequality. For example,
what would be the optimal constant $C$ in the inequality (6.6)?
\endroster
\enddefinition

\medskip

\head{\bf Appendix 1: Smoothness implies Hamiltonian continuity}
\endhead
\smallskip

We first recall the precise definition of smooth Hamiltonian
paths.

\definition{Definition A.1}
(i) A $C^\infty$-diffeomorphism $\phi$ of $(M,\omega)$ is a {\it
Hamiltonian diffeomorphism} if $\phi = \phi_H^1$ is the time-one map of the Hamilton equation
$$ \dot x = X_H(t,x) $$
for a $C^\infty$ function $H: \R \times M \to \R$ such that
$$ H(t+1,x) = H(t,x) $$
for all $(t,x) \in \R \times M$.
We denote by $Ham(M,\omega)$ the
set of Hamiltonian diffeomorphisms with the $C^\infty$-topology
induced by the inclusion
$$
Ham(M,\omega) \subset Symp_0(M,\omega),
$$
where $Symp_0(M,\omega)$ carries the $C^\infty$-topology.
\smallskip

\n (ii) A {\it (smooth) Hamiltonian path} $\lambda: [0,1] \to
Ham(M,\omega)$ is a smooth map
$$
\Lambda: [0,1] \times M \to M
$$
such that \roster
\item its derivative $\dot\lambda(t) = \frac{\part \lambda}{\part t}\circ
(\lambda(t))^{-1}$ is Hamiltonian, i.e., the one form
$\dot\lambda(t) \, \rfloor \omega$ is exact for all $t \in [0,1]$.
We call a function $H: \R \times M \to \R$ a
{\it generating Hamiltonian} of $\lambda$ if it satisfies
$$
\lambda(t) = \phi_H^t \circ \lambda(0), \quad \text{or equivalently,} \quad
dH_t = \dot\lambda(t) \, \rfloor \omega.
$$
\item
$\lambda(0):=\Lambda(0, \cdot): M \to M$ is a Hamiltonian
diffeomorphism, and therefore $\lambda(t) = \Lambda(t, \cdot)$ is for all
$t\in [0,1]$.
\endroster
We denote by $\PP^{ham}(Symp(M,\omega))$ the set of Hamiltonian
paths $\lambda: [0,1] \to Ham(M,\omega)$, and by
$\PP^{ham}(Symp(M,\omega),id)$ the set of such $\lambda$ with
$\lambda(0) = id$. We provide the obvious topology on
$\PP^{ham}(Symp(M,\omega))$ and $\PP^{ham}(Symp(M,\omega),id)$
induced by the $C^\infty$-topology of the space $C^\infty ([0,1] \times M, M)$ of
the corresponding maps $\Lambda$ above. We call this the $C^\infty$-topology of
$\PP^{ham}(Symp(M,\omega))$ and $\PP^{ham}(Symp(M,\omega),id)$.
\enddefinition

Note that if $\phi = \phi_H^1$ is a Hamiltonian diffeomorphism (in
the sense of definition A.1.(i)), then $t \mapsto \lambda(t) =
\phi_H^t$ is a smooth Hamiltonian path (in the sense of definition
A.1.(ii)) with $\lambda(0) = id$ and $\lambda(1) = \phi$. So each
$\phi \in Ham(M,\omega)$ can be connected to the identity by a
smooth Hamiltonian path as in A.1.(ii). In particular,
$Ham(M,\omega)$ is the image of the evaluation map $ev_1$ (1.5).
We also note that by Proposition 3.4, each smooth path
$\lambda:[0,1] \to Symp(M,\omega)$ that has its image contained in
$Ham(M,\omega)$ is a smooth Hamiltonian path in the sense of
Definition A.1 (ii).

 In this appendix, we give the proof of the following
basic lemma and prove that any smooth path in $Ham(M,\omega)$ is
Hamiltonian continuous. By abuse of notation, we will just denote
a smooth Hamiltonian path by
$$
\lambda: I \to Ham(M,\omega),
$$
or more generally, a smooth Hamiltonian map from a simplex $\Delta$
by
$$
\lambda : \Delta \to Ham(M,\omega).
$$

\proclaim{Lemma A.2} For any Hamiltonian path $\lambda: I \to
Ham(M,\omega)$ defined on an interval $I = [a,b]$ such that $\lambda$ is
flat near $a$, i.e., there exists $a' > a$ with
$$
\lambda(s) \equiv \lambda(a) \tag A.1
$$
for all $ a \leq s \leq a' \leq b$, we can find a
smooth map
$$
\Lambda: I \times [0,1] \times M \to M
$$
such that the following hold: \roster
\item For each $s \in I$ and $t \in [0,1]$, $\Lambda_{(s,t)} \in
Ham(M,\omega)$, where we denote
$$ \Lambda_{(s,t)}(x): = \Lambda(s,t,x). $$
\item For each $s \in I$, the path $\lambda^s: [0,1] \to
Ham(M,\omega)$ is a Hamiltonian path with $\lambda^s(0) = id$ and $\lambda^s(1) = \lambda(s)$,
 which is flat near $0$,
where we denote
$$\lambda^s (t) := \Lambda_{(s,t)}. $$
\endroster
Furthermore, a similar statement holds for a map $ \Delta \to
Ham(M,\omega)$ where $\Delta$ is a $k$-simplex: in this case (A.1)
is replaced by the condition that $\lambda$ is flat near the
vertex $0 \in \Delta$.
\endproclaim
\demo{Proof} We may assume $I = [0,1]$.
Let $K: I \times M \to \R$ be the (not necessarily normalized)
Hamiltonian generating $\lambda$ such that
$$
\lambda(s) = \phi_K^s \circ \lambda(0), \quad s\in [0,1] \tag A.2
$$
and
$$
K(s, \cdot) \equiv 0 \quad \text{for all} \quad 0 \leq s \leq a'. \tag A.3
$$
(A.3) is possible because of the assumption (A.1). Next we fix a
Hamiltonian $H^0 : [0,1] \times M \to \R$ with $H^0 \mapsto
\lambda(0)$. After reparameterization, we may assume that
$$
H^0 \equiv 0 \quad \text{near} \quad t = 0, \, 1. \tag A.4
$$
Now for each $s \in [0, 1]$, we define $H^s: [0,1] \times M \to
\R$ by the formula
$$
H^s(t,x) = \cases \frac{1}{1 - s}H^0\left(\frac{1}{1-s} t, x\right)& \quad
\text{for } 0 \leq t < 1-s, \\
K(t-(1-s), x) & \quad \text{for }  1 - s \leq t \leq 1 .
\endcases
\tag A.5
$$
Obviously $H: I \times [0,1] \times M \to \R$ is smooth due to the above
flatness conditions (A.3) and (A.4) and satisfies
$$ \phi_{H^s}^1 = \lambda(s). $$
We then define $\Lambda$ by $\Lambda(s,t) = \phi_{H^s}^t$. It
follows from the construction that $\Lambda$ satisfies all the
properties in (1) and (2). The last statement can be proven by a similar argument by
considering the retraction of the $k$-simplex $\Delta$ to its vertex $0$. \qed\enddemo

Remark that if $\lambda$ is flat also near $t=1$,
then we can assume that $\lambda^s$ is flat near $t=1$ for all $s \in I$.
The proof goes through the same way.

\proclaim{Corollary A.3} Any smooth Hamiltonian path $\lambda :
[0,1] \to Ham(M,\omega)$ is Hamiltonian continuous.
\endproclaim
\demo{Proof} Let $\lambda = \phi_H : [0,1] \longrightarrow Ham
(M,\omega)$ be a smooth Hamiltonian path (in the sense of
Definition A.1.(ii)). Here we assume without loss of generalities
that $\lambda(0) = id$. We have to show that $\lambda$ is
continuous with respect to the Hamiltonian topology on
$Ham(M,\omega)$, i.e., as a map $\lambda : [0,1] \longrightarrow
\HH am (M,\omega)$. Note that $\lambda$ factors through
$$ [0,1] \to \PP_s^{ham}(Symp(M,\omega),id) \to \HH am(M,\omega),
\quad s \mapsto \phi_{H^s} \mapsto \phi_{H^s}^1 = \phi_H^s
=\lambda(s),
$$ where the second map is the evaluation map. By definition of
the Hamiltonian topology on $Ham(M,\omega)$, it suffices to prove
that the first map is continuous. The topology on
$\PP_s^{ham}(Symp(M,\omega),id)$ is by Proposition 3.10 equivalent
to the metric topology induced by $d_{ham}$. So we only have to
show that the map $s \mapsto \phi_{H^s}$ is continuous with
respect to the standard metric on $[0,1]$ and $d_{ham}$ on
$\PP_s^{ham}(Symp(M,\omega),id)$.

Let $H^s$ be the Hamiltonian and $\Lambda$ be the smooth map
constructed in the proof of Lemma A.2. By definition
$$ d_{ham}\left(\phi_{H^s},\phi_{H^{s'}}\right) = \| H^{s} - H^{s'} \|
+ \overline d \left(\phi_{H^s},\phi_{H^{s'}}\right) . \tag A.6$$
If we define the smooth reparameterization functions $\zeta_1$,
$\zeta_2 : [0,1] \to [0,1]$, $\zeta_1 (t) = st$, $\zeta_2 (t) = s'
t$, then $\| \zeta_1 - \zeta_2 \|_{ham} = 2 | s - s' |$. Hence by
Lemma 3.20, the first term in (A.6) is less than $ 2 C | s - s'
|$, where $C$ is the constant given in (3.20) in Lemma 3.20. For
the second term in (A.6), first note that $\Lambda$ is Lipschitz
continuous since it is smooth and compactly supported. Therefore,
$$
d_{C^0} \left( \phi_{H^s},\phi_{H^{s'}} \right) = \max_{(t,x)}
d \Big( \Lambda (s,t,x), \Lambda(s',t,x) \Big) < L | s - s' |,
$$
where $L$ is a Lipschitz constant for $\Lambda$.
Since $s \mapsto (\lambda(s))^{-1}$ is also a smooth Hamiltonian path,
we can use Lemma A.2 to construct a corresponding map $\Lambda' (s,t) = (\phi_{H^s}^t)^{-1}$,
and then apply the same argument to obtain
$$
d_{C^0} \left( (\phi_{H^s})^{-1},(\phi_{H^{s'}})^{-1} \right) < L' | s - s'| , $$
where $L'$ is another Lipschitz constant.
That shows that the second term in (A.6) is less than $\max (L,L') | s - s' |$.
Altogether, with $c = \max (2C,L,L')$, we have
$$
d_{ham}\left(\phi_{H^s},\phi_{H^{s'}}\right) = \| H^{s} - H^{s'}
\| + \overline d (\phi_{H^s},\phi_{H^{s'}}) < c | s - s' | , $$
which completes the proof. \qed \enddemo

\head{\bf Appendix 2: The $L^{(1,\infty)}$-Approximation Lemma}
\endhead
\smallskip

In this appendix, we give the proof of the $L^{(1,\infty)}$-Approximation
Lemma which is a slight variation of [Lemma 5.2, Oh3].

\proclaim{Lemma A.4 ($L^{(1,\infty)}$-Approximation Lemma)} Let $H
: [0,1] \times M \to \R$ be a given smooth Hamiltonian and $\phi =
\phi_H^1$ be its time-one map. Then we can reparameterize
$\phi_H^t$ in time so that the Hamiltonian $H'$ generating the
reparameterized path satisfies the following properties: \roster
\item $\phi_{H'}^1 = \phi_H^1$,
\smallskip

\item $H' \equiv 0$ near $t = 0, \, 1$, and in particular $H'$ can be
extended to be time-periodic on $\R \times M$,
\smallskip

\item the norm $\| \overline H\#H' \|$ can be made as small as we
want, and
\smallskip

\item for the Hamiltonians $H', \, H''$ generating any two such
reparameterizations of $\phi_H^t$, there is a canonical one-one
correspondence between $\text{Per}(H')$ and $\text{Per}(H'')$,
and $\text{Crit }\AA_{H'}$ and $\text{Crit }\AA_{H''}$ with their
actions fixed .
\endroster
Furthermore this reparameterization is canonical in the sense
that the ``smallness'' in (3) can be chosen uniformly over $H$
depending only on the $C^0$-norm and the modulus of continuity of
$H$. In particular, this approximation can be done with respect to
the Hamiltonian topology. Moreover, the closeness in the
Hamiltonian topology can be made as small as we want
independent of $H$ (only the time for which the reparameterized
Hamiltonian is flat depends on $H$).
\endproclaim
\demo{Proof}
We first reparameterize $\phi^t_H$ in the following way: We choose
a smooth function $\zeta : [0,1] \longrightarrow [0,1]$ such that for $\e > 0$
$$
\zeta (t) = {\cases
             0  & \text {for }\, 0 \leq t \leq \epsilon \\
             1  &  \text {for }\, 1 - \epsilon \leq t \leq 1
             \endcases}
$$
and
$$
\zeta'(t) \geq 0 \quad \hbox {for all} \quad t \in [0,1],
$$
and consider the isotopy
$$
\psi^t : = \phi_H^{\zeta (t)} .
$$
It is easy to check that the Hamiltonian generating the isotopy
$\{ \psi^t \}_{0 \leq t \leq 1} $ is $H' = \{ H'_t \}_{0 \leq t
\leq 1}$ with $H'_t = \zeta '(t) H_{\zeta (t)}$.  By definition,
it follows that $H'$ satisfies (1) and (2).
As always we assume that $H$ is normalized, and then so is $H'$.
In particular, $\int^1_0 \max (H' - H) dt \ge 0 $.
For (3), we compute
$$
\align
    0 \le \int^1_0 \max_x & (H' - H) dt = \int^1_0 \max_x (\zeta ' (t)
H_{\zeta(t)} - H_t ) dt \\
       &  \le \int^1_0
\max_x \Big( \zeta ' (t) (H_{\zeta (t)} - H_t )\Big) dt +
                  \int^1_0 \max_x \Big(( \zeta ' (t) - 1) H_t\Big)dt .
\endalign
$$
For the first term,
$$
\align \int^1_0 \max_x & \Big(\zeta '(t) (H_{\zeta (t)} -
H_t)\Big)dt =\int^1_0 \zeta '(t) \max_x (H_{\zeta (t)} - H_t)dt\\
& \leq \int^1_0 \zeta '(t) \max_{x,t} | H_{\zeta(t)} - H_t | dt =
\max_{x,t} |H_{\zeta (t)} (x) - H_t (x) | \leq L \cdot \| \zeta  - id \|_{C^0}
\endalign
$$
which can be made arbitrarily small by choosing $\zeta$ so that
$\| \zeta  - id \|_{C^0}$ becomes sufficiently small. Here $L$ is
a Lipschitz constant for $H$ in the time variable $t$ (it exists
and is finite since $H$ is smooth and supported on the compact set
$[0,1] \times M$). We refer to this constant as the modulus of
continuity. For the second term,
$$ \int^1_0 \max_x \Big(( \zeta ' (t) - 1) H_t\Big) dt
\leq \int^1_0 | \zeta '(t) -1 | dt \cdot \max_{x,t}|H(x,t)| = \|
H\|_{C^0} \int^1_0 | \zeta ' (t) - 1 |dt. $$ Again by
appropriately choosing $\zeta$ (which can be done consistently
with the choice above), we can make
$$
\int^1_0 | \zeta '(t) - 1 | dt
$$
as small as we want.  Combining these two, we have verified
$\int^1_0 \max (H'-H)\, dt$ can be made as small as we want
by making the hamiltonian norm
$$
\|\zeta - id\|_{ham} = \| \zeta  - id \|_{C^0} + \| \zeta ' - 1 \|_{L^1}
$$
small. This can always be done by choosing $\e$ sufficiently
small. Similar consideration applies to $\int^1_0 - \min (H' - H)
\, dt$ and hence we have finished the proof of (3).

The statement (4) follows from simple comparison of the
corresponding actions of periodic orbits. The statements in the
last paragraph follow from the construction. For the
$C^0$-closeness, note that similarly to the proof above, by
continuity of the path $t \mapsto \phi_H^t$, the distance
$\overline{d} \left( \phi_{H^\zeta}, \phi_H \right)$ can be made
arbitrarily small by choosing $\zeta$ so that $\| \zeta - id
\|_{C^0}$ becomes small. This finishes the proof. \qed\enddemo

We would like to point out that the above modification does {\it
not} approximate in the $L^\infty$-topology on $[0,1] \times M$
because the derivative of the cut-off function $\zeta$ could blow
up in the above approximation. In fact it is easy to see that such
an approximation can be done for a given Hamiltonian function $H$
in the $L^\infty$-norm if and only if $H_0 \equiv H_1 \equiv$
constant. The proof is essentially the same as above.

\demo{Proof of Lemma 3.20} Replace $\zeta$ by $\zeta_1$ and $id$
by $\zeta_2$ in the proof of the $L^{(1,\infty)}$-Approximation
Lemma. \qed
\enddemo

\head {\bf References}
\endhead

\enddocument